\documentclass[12pt]{article}
\usepackage{pb-diagram}
\usepackage{lamsarrow}
\usepackage{pb-lams}
\usepackage{amsmath}
\usepackage{amssymb}
\newtheorem{theo}{Theorem}
\newtheorem{prop}{Proposition}
\newtheorem{lemm}{Lemma}
\newtheorem{coro}{Corollary}
\newtheorem{rema}{Remark}

\newtheorem{conj}{Conjecture}

\newcommand{\cqfd}
{%
\mbox{}%
\nolinebreak%
\hfill%
\rule{2mm}{2mm}%
\medbreak%
\par%
}
\newfont{\gothic}{eufb10}
\date{\empty}
\begin{document}
\title{Green's canonical syzygy conjecture for generic curves of odd genus}
\author{Claire Voisin\\ Institut de math{\'e}matiques de Jussieu, CNRS,UMR
7586} \maketitle \setcounter{section}{-1}
\section{Introduction}
For $X$ a projective variety, $L$ a line bundle on $X$, and
$\mathcal{F}$ a coherent sheaf on $X$, denote, following
\cite{Gre}, by $K_{p,q}(X,L,\mathcal{F})$ the cohomology at the
middle of the exact sequence
$$\bigwedge^{p+1}H^0(X,L)\otimes H^0(X,\mathcal{F}((q-1)L))\rightarrow
\bigwedge^{p}H^0(X,L)\otimes H^0(X,\mathcal{F}(qL))$$
$$\rightarrow
\bigwedge^{p-1}H^0(X,L)\otimes H^0(X,\mathcal{F}((q+1)L)),$$ where
the maps are  Koszul differentials. For
$\mathcal{F}=\mathcal{O}_X$, use the notation $K_{p,q}(X,L)$.
Green's conjecture on syzygies of canonical curves (see
\cite{Gre}) relates the Koszul cohomology groups
$$K_{p,1}(C,K_C),$$
for $C$ a smooth projective curve, to the Clifford index of the
curve :
$$Cliff(C):=Min_{D}\{deg\,D-2r\},$$
where $D$ runs through the set of divisors $D$ on $C$ satisfying :
$$r+1:=h^0(D)\geq2,\,h^1(D)\geq2.$$
\begin{conj}\label{greenconj} (Green)
$$K_{l,1}(C,K_C)=0,\,\forall l\geq p\Leftrightarrow Cliff(C)>g-p-2.$$
\end{conj}
The direction $ \Rightarrow$ is proved by Green and Lazarsfeld in
the appendix to \cite{Gre}.
 The case $p=g-2$ of the conjecture is equivalent to Noether's
 theorem, and
 the case $p=g-3$ to Petri's theorem (see \cite{georoma}). The
 case $p=g-4$ has been proved in any genus by Schreyer \cite{Sch2} and by the
 author \cite{Vo2} for $g>10$.

 More recently, the conjecture has been studied in \cite{teduke},
 \cite{Vo}, for generic curves of fixed gonality.
 Teixidor proves the following
 \begin{theo}  (\cite{teduke}) Green's conjecture is true for generic curves of
 genus $g$ and fixed gonality $\gamma$, in the range
 $$\gamma\leq \frac{g+7}{3}. $$
 \end{theo}
 Note that Brill-Noether theory says that the gonality $\gamma$
 always satisfies the inequality
 $$\gamma\leq [\frac{g+3}{2}],$$
 with equality for the generic curve.
 We proved  the following
 \begin{theo} (\cite{Vo})
 Green's conjecture is true for generic curves of
 genus $g$ and fixed gonality $\gamma$, in the range
 $$\gamma\geq \frac{g}{3}, $$
 except possibly for the generic curves of odd genus $g=2k+1$, whose gonality is
 $k+2$.
 \end{theo}
 So, for generic curves of
 fixed gonality,  the only remaining case is that of generic curves of odd genus $g=2k+1$.
 Green's conjecture together with Brill-Noether theory predicts
 that
 \begin{eqnarray}\label{stat}
 K_{k,1}(C,K_C)=0.
 \end{eqnarray}
 This is the main result proved in this paper. We give the precise
 statement below; it gives slightly more, since it proves the
 vanishing (\ref{stat}) for some explicit curves which we know to be
 generic in the Brill-Noether sense. Applications of this result
 to the gonality conjecture for generic curves
 of even genus can be found in \cite{AV}, \cite{aprodu}.

 Note that this last case was especially challenging, first of all
 because, as noticed in
 \cite{HR}, the locus of jumping syzygies, i.e. the locus where
$K_{k,1}(C,K_C)\not=0$ is of codimension $1$ in $\mathcal{M}_g$ in
this case, and in fact has a natural structure of determinantal
hypersurface,
 and also because of the following important result of Hirschowitz
 and
Ramanan :
\begin{theo}(\cite{HR}) If the Green conjecture is true
for generic curves of genus $2k+1$, then  the locus of jumping
syzygies in $\mathcal{M}_{2k+1}$ is equal set theoretically to the
$k+1$-gonal divisor, which is also the locus where the Clifford
index is one less the generic Clifford index.
\end{theo}
Combined with the generic Green conjecture for genus
$2k+1$-curves, this  provides a strong evidence for   conjecture
\ref{greenconj}.

Coming back to our result, the curves we consider are the
following : we consider a smooth projective $K3$ surface $S$, such
that $Pic(S)$ is isomorphic to $\mathbb{Z}^2$, and is freely
generated by $L$ and $\mathcal{O}_S(\Delta)$, where $\Delta$ is a
smooth rational curve such that $deg\,L_{\mid \Delta}=2$, and $L$
is a very ample line bundle
 with $L^2=2g-2, \,g=2k+1$. By the hyperplane section theorem (see
 \cite{Gre}), we have
 $$K_{k,1}(S,L)\cong K_{k,1}(C,L),\,\forall C\in\mid L\mid.$$
 As we shall see in the next section, curves in $\mid L\mid$ have
 the generic Clifford index. Hence we expect from Green's conjecture that
 $$K_{k,1}(C,L)=0=K_{k,1}(S,L).$$
 Our theorem says indeed :
 \begin{theo}\label{ourthm} The $K3$ surface $S$ being as above, we have
 $$K_{k,1}(S,L)=0.$$
 \end{theo}
 In the first section, we show how to adapt the arguments of
 \cite{Vo}, to the line bundle $L+\Delta$ on $S$, in order to show
 that
 $$K_{k+1,1}(S,L+\Delta)=0.$$
 Note that the proof of \cite{Vo} worked under the assumption
 $Pic(S)=\mathbb{Z}$, which is why a few supplementary arguments
 are needed.

 In the second section, we show how to deduce from this the
 vanishing $K_{k,1}(S,L)=0$. The last section is devoted to the
 proof of the crucial proposition \ref{injective} used in the
 proof of Theorem \ref{ourthm}.

 \vspace{0,5cm}

 {\bf Acknowledgements:} I am very indebted to the referee for
 pointing out a mistake in my original proof of Proposition
 \ref{injective} and for helpful comments.

\section{The case of curves of even genus on a $K3$ surface with a node}
Let $S$ be a $K3$ surface, whose Picard group is freely generated
by a very ample line bundle $L$, such that
$$L^2=2g-2,\,g=2k+1,$$
and $\mathcal{O}_S(\Delta)$, where $\Delta$ is a rational curve
such that
$$deg\,L_{\mid\Delta}=2.$$
Let $L'=L(\Delta)$ ; smooth curves in $\mid L'\mid$ do not meet
$\Delta$ and are of genus $2k+2=2(k+1)$. Contracting $\Delta$ to a
node, the line bundle $L'$ descends, and we are essentially in the
situation considered in \cite{Vo}. (Note however the change of
notations  from $k$ to $k+1$.)

We first apply Lazarsfeld's argument in \cite{La} to show :
\begin{prop}\label{BN} Smooth curves $C$ in $\mid L\mid $ or in
$\mid L'\mid$ are generic in the Brill-Noether sense, i.e., do not
have a $g_d^r$ when the Brill-Noether number $\rho(g(C),d,r)$ is
negative. In particular, their Clifford index is the generic one.
\end{prop}
{\bf Proof.} It follows from \cite{La} that if $C\subset S$ is a
smooth curve in a linear system $\mid M \mid$, and $D$ is a
$g_d^r$ on $C$ with $\rho(g(C),d,r)<0$, there exists a line bundle
$H$ on $S$ with
$$h^0(H)\geq2,\,h^0(M-H)\geq2.$$
Apply this to $M=L$ or $M=L'$. Writing $H=\alpha L+\beta\Delta$,
the condition $h^0(H)\geq2$ implies that $\alpha\geq1$. Similarly,
the condition $h^0(M-H)\geq2$ implies that $1-\alpha\geq1$,
whether $M=L$ or $M=L'$. This is a contradiction. \cqfd It is now
expected from Green's conjecture \ref{greenconj} that
$$K_{k+1,1}(S,L')=0,\,K_{k,1}(S,L)=0.$$
In \cite{Vo}, we proved the vanishing
$$K_{k+1,1}(S,L')=0,$$
for a line bundle $L'$ on $S$ with
$${L'}^2=2g'-2,\,g'=2k+2,$$
under the assumption that $L'$ generates $Pic(S)$. Our first goal
is to extend this result in our situation.
\begin{theo} \label{evennode}For $S,\,L'=L(\Delta)$ as above, we have
\begin{eqnarray}\label{vanth}K_{k+1,1}(S,L')=0.\end{eqnarray}
\end{theo}
The proof of this theorem occupies the rest of this section. Let
$C'\in\mid L'\mid$ be smooth; by Brill-Noether theory, there is a
smooth $g_{k+2}^1$, say $D$, on $C'$. By proposition \ref{BN},
both $D$ and $K_{C'}-D$ are generated by sections. Consider the
Lazarsfeld bundle
\begin{eqnarray}\label{Lmuk}
E=F^*=F\otimes L',
\end{eqnarray}
where $F$ is the rank $2$ vector bundle fitting in the exact
sequence
\begin{eqnarray}\label{exE}
0\rightarrow F\rightarrow H^0(C',D)\otimes\mathcal{O}_S\rightarrow
D\rightarrow0.
\end{eqnarray}
Here the last map is the evaluation map along $C'$.
 One can show that $E$ does not depend on the curve $C'$, and  neither on
 $D'$. The bundle $E$ has $det\,E\cong L'$, and $h^0(E)=k+3$.
 The following key point, which was used constantly throughout the
 proof of
 \cite{Vo}, remains true in our situation :
 \begin{prop}\label{detgrass} The determinant map
 $$d:\bigwedge^2H^0(S,E)\rightarrow H^0(S,L')$$
 does not vanish on decomposable elements.
 \end{prop}
{\bf Proof.} Indeed assume $s,\,s'\in H^0(S,E)$ are not
 proportional but satisfy
 $d(s\wedge s')=0$. Then $s,\,s'$ generate a sub-line bundle
 of $E$, say $H$, which we may assume saturated, and which
 satisfies
 $$h^0(H)\geq2.$$
 Hence there is an exact sequence
 $$0\rightarrow H\rightarrow E\rightarrow H'\rightarrow
 T\rightarrow0,$$
 where $H'$ is a line bundle such that $H+H'=det\,E=L'$, and $T$ is
 torsion supported on points of $S$.
 Since $E$ is generated by sections, $H'$ is generated by
 sections away from the support of $T$. On the other hand
 $H'$ is not trivial, since $H^0(S,E^*)=0$. So $h^0(H')\geq2$.
 But this contradicts the fact we already mentioned, that we cannot write $L'$ as the sum of
 two
 line bundles admitting at least two sections.
\cqfd We now recall the main points of the proof of the vanishing
(\ref{vanth}) given in \cite{Vo}, in order to make clear what has
to be added in our situation. We warn again the reader that the
notation of \cite{Vo} has been shifted (the integer $k$ there
becomes $k+1$ here).

{\it First step.}  Let $S^{[k+2]}_{curv}$ be the open subset of
the Hilbert scheme of $S$ parametrizing curvilinear,  degree
$k+2$, $0$-dimensional subschemes of $S$. Let
$$I_{k+2}\stackrel{\pi_{k+2}}{\rightarrow} S^{[k+2]}_{curv},\,
I_{k+2}\subset S\times S^{[k+1]}_{curv}$$ be the incidence scheme.
We established the following isomorphism :
\begin{eqnarray}
\label{step1} K_{k+1,1}(S,L')\cong
H^0(I_{k+2},\pi_{k+2}^*L'_{k+2})/\pi_{k+2}^*
H^0(S^{[k+2]}_{curv},L'_{k+2}),
\end{eqnarray}
where the line bundle $L'_{k+2}$ is the determinant of the vector
bundle ${\mathcal E}_{L'}$ of rank $k+2$ on $S^{[k+2]}_{curv}$,
defined as
$${\mathcal E}_{L'}=R^0\pi_{k+2*}(pr_1^*L').$$
From this we deduced the following criterion :

\begin{lemm} \label{istep1} The vanishing $K_{k+1,1}(S,L')=0$ holds if there exists a
reduced  scheme $Z$, and a morphism
$$j:Z\rightarrow S^{[k+2]}_{curv}$$
such that, denoting
$$\widetilde Z \stackrel {\tilde j}{\rightarrow}
I_{k+2}$$ the fibered product
$$\widetilde Z=Z\times_{S^{[k+2]}_{curv}}I_{k+2}$$
we have :
\begin{enumerate}
\item   The map
\begin{eqnarray} \tilde j^*:H^0(I_{k+2},\pi_{k+2}^*L'_{k+2})
\rightarrow H^0(\widetilde Z,\tilde j^*\pi_{k+2}^*L'_{k+2})
\label{j*}
\end{eqnarray}
is injective. \item \label{istep2} Denoting by $\pi:\widetilde
Z\rightarrow Z$ the first projection, the map
\begin{eqnarray}\label{pi*}
\pi^*:H^0(Z,j^*L'_{k+2})\rightarrow H^0(\widetilde Z, \tilde
j^*\pi_{k+2}^*L'_{k+2})
\end{eqnarray}
is surjective.
\end{enumerate}
\end{lemm}
{\it Second step.} The construction of $Z$ is as follows : we
start with the vector bundle $E$ of  (\ref{Lmuk}), (\ref{exE}). It
has $c_2(E)=k+2$. Denote by
$${\mathbb P}(H^0(E))_{curv}\subset {\mathbb P}(H^0(S,E))$$
the open set parametrizing sections $\sigma\in H^0(S,E)$ whose
$0$-scheme $z_\sigma$ is $0$-dimensional and curvilinear.
 There is a natural morphism
$${\mathbb P}(H^0(E))_{curv}\rightarrow S^{[k+2]}_{curv}$$
$$\sigma\mapsto z_\sigma.$$
Let $W={\mathbb
P}(H^0(E))_{curv}\times_{S^{[k+2]}_{curv}}I_{k+2}$. This is a
degree $k+2$ cover of ${\mathbb P}(H^0(E))_{curv}$. It admits a
natural morphism, say $f$ to $I_{k+2}$.
 We use now the morphism
$$\tau_{k+2} :I_{k+2}\rightarrow S\times S^{[k+1]}_{curv},$$
which sends a point $(x,z),\{x\}\subset  z$ of $I_{k+2}$ to the
residual scheme of $x$ in $z$, which is curvilinear of length
$k+1$, since $z$ is curvilinear. Let
$$\psi:W\rightarrow S^{[k+1]}_{curv}$$
be the composed map $\psi= \tau_{k+2} \circ f$. Finally we
construct the sum map:
$$j:Z:=(\widetilde{S\times W})_0\rightarrow S^{[k+2]}_{curv},$$
$$(x,w)\mapsto x\cup \psi(w),$$
where the    $\,\,\,\widetilde{}\,\,\,$   stands here to mean
``blowup along the incidence subscheme in order to make the scheme
structure on the union $x\cup\psi(w)$ well defined'', and the
subscript $0$ means, ``taking an open set in order to make sure
that this scheme structure is curvilinear''.

{\it Third step .}  The injectivity of the map (\ref{j*}) in Lemma
\ref{istep1}  is easily reduced to the injectivity of the
restriction map
$$\psi^*:H^0(S^{[k+1]}_{curv},L'_{k+1})\cong \bigwedge^{k+1}
H^0(S,L')\rightarrow H^0(W,\psi^*L'_{k+1}).$$ Now if
$\beta:W\rightarrow {\mathbb P}(H^0(E))_{curv}$ is the natural
surjective map, we showed that
$$\psi^*L'_{k+1}\cong\beta^*{\mathcal O}_{ {\mathbb P}(H^0(E))_{curv}}(k+1)$$
and that the map above is the composition of $\beta^*$ and of an
isomorphism
$$\bigwedge^{k+1}
H^0(S,L')\cong H^0( {\mathbb P}(H^0(E))_{curv},{\mathcal O}(k+1))
\cong S^{k+1}H^0(S,E)^*.$$ The construction of this isomorphism
uses only the proposition \ref{detgrass} which remains true in our
situation. Hence this  step works as in \cite{Vo}.

{\it Fourth step.}
 In \cite{Vo}, we reduced easily the proof of the surjectivity of the map
(\ref{pi*}) in Lemma \ref{istep1}, to the proof of the following :
let
$$\widetilde W=W\times_{S^{[k+1]}} I_{k+1},$$
and denote by $\gamma:\widetilde W\rightarrow W$ the natural map.
\begin{prop}
The map
\begin{eqnarray}\label{gamma*}
\gamma^*:H^0(W,\psi^*L_{k+1})\rightarrow H^0(\widetilde
W,\gamma^*\psi^*L_{k+1})
\end{eqnarray}
is surjective.
\end{prop}
Using the fact that
$$\psi^*L_{k+1}=r^*{\mathcal O}_{{\mathbb P}(H^0(E))_{curv}}(k+1),$$
where
$$r=\beta\circ \gamma:\widetilde W\rightarrow
{\mathbb P}(H^0(E))_{curv},$$ this proposition is a consequence of
the following :
\begin{prop}\label{tresdur}
The map
\begin{eqnarray}
r^*:H^0({\mathcal O}_{{\mathbb P}(H^0(E))_{curv}}(k+1))=
S^{k+1}H^0(S,E)^*\rightarrow H^0(\widetilde W,r^*{\mathcal
O}_{{\mathbb P}(H^0(E))_{curv}}(k+1))
\end{eqnarray}
is surjective.
\end{prop}
This is in the proof of this proposition that we shall see a
difference between the case considered in \cite{Vo} and the
present case. Indeed, let us introduce as in \cite{Vo}, the
codimension $4$ subscheme
\begin{eqnarray}\label{W'}
W' =\{(z,\sigma)\in \widetilde{S\times S}\times {\mathbb
P}(H^0(S,E)),\,\,\sigma_{\mid z}=0 \},
\end{eqnarray}
where $\widetilde{S\times S}$ is the blowup of $S\times S$ along
its diagonal, hence parametrizes ordered length $2$ subschemes of
$S\times S$.

In \cite{Vo}, we used the fact that $\widetilde{W}$ can be seen as
a large (i.e. the complementary set has codimension $\geq2$)
Zariski open set in $W'$, and the fact (which is inaccurately not
mentioned explicitly) that $W'$ is normal (in fact it is smooth
for $k+1>3$, see below) to conclude that
$$H^0(\widetilde W,r^*{\mathcal O}_{{\mathbb P}(H^0(E))_{curv}}(k+1))=
H^0(W',pr_2^*{\mathcal O}_{{\mathbb P}(H^0(S,E))}(k+1)).$$

Here we cannot do that because it is not true anymore that $W'$ is
normal, nor that $\widetilde W$ is large in $W'$. In fact $W'$ is
not irreducible. Indeed, consider the rational curve
$\Delta\subset S$. The exact sequence (\ref{exE}) together with
the fact that $L'_{\mid \Delta}$ is trivial, shows that $E_{\mid
\Delta}$ is trivial and that the restriction map
$$H^0(S,E)\rightarrow H^0(\Delta,E_{\mid \Delta})$$
is surjective, the right hand side being of rank $2$. So
$H^0(S,E(-\Delta))$ is of codimension $2$ in $H^0(S,E)$, so that
$W'$ has one component isomorphic to $\Delta\times\Delta\times
{\mathbb P}(H^0(S,E(-\Delta))$.

However, what remains true in our situation is the following
\begin{lemm}\label{smooth} Assume $k+1>3$ (Theorem
\ref{evennode} is already known for $k+1\leq3$).  Let
$$U:={\mathbb P}(H^0(S,E))-{\mathbb P}(H^0(S,E(-\Delta))).$$
Then $W'_U:=W'\cap (\widetilde{S\times S}\times U)$ is smooth and
$\widetilde{W}$ is a large open set in it.
\end{lemm}
{\bf Proof.} If a section of $E$ vanishes at a point of $\Delta$,
then it vanishes along $\Delta$, hence its $0$-locus is not
$0$-dimensional. So $\widetilde W\subset W'$. The proof that it is
a large open set is easy.
 To prove that $W'_U$ is smooth, it suffices to show the following :
 \newline\newline
 {\it For any $z\in (S-\Delta)^{[2]}$, the restriction map
\begin{eqnarray}\label{restr}
H^0(S,E)\rightarrow H^0(E_{\mid z})
\end{eqnarray}
is surjective.}
\newline\newline
Choose a smooth curve $C'\in \mid L'\mid$ containing $z$. It
exists because $z$ does not meet $\Delta$. There is an exact
sequence
$$0\rightarrow D\rightarrow E_{\mid C'}\rightarrow
K_{C'}-D\rightarrow 0,$$ where $D$ is a divisor of degree $k+2$ on
$C'$, with $h^0(D)=2$. Furthermore the map $H^0(S,E)\rightarrow
H^0(C',E_{\mid C'})$ is surjective and there is an exact sequence
$$0\rightarrow H^0(C',D)\rightarrow H^0(C',E_{\mid C'})
\rightarrow H^0(C',K_{C'}-D)\rightarrow0.$$ Now, by proposition
\ref{BN}, the curve $C'$ is generic in the Brill-Noether sense.
Hence, since $k+1>3$, it does not possesses a $g_{k+4}^2$. Hence
the map
$$H^0(C',K_{C'}-D)\rightarrow H^0((K_{C'}-D)_{\mid z})$$
is surjective. So the map (\ref{restr}) has at least rank $3$
since $\mid D\mid$ has no base point so that the restriction map
$H^0(C',D)\rightarrow H^0(D_{\mid z})$ has at least rank $1$, and
our statement will be proved if we can furthermore choose $C'$ and
$D$ so that the restriction map
$$H^0(C',D)\rightarrow H^0(D_{\mid z})$$ is injective.
Take now two sections $s,\,s'$ of $E$ such that $d(s\wedge s')$
vanishes on $z$, but $s,s'$ have independent restrictions in
$H^0(E_{\mid z})$. It is easily shown to exist once we know that
the restriction map (\ref{restr}) has rank at least $3$. Let $C'$
be defined by $d(s\wedge s')$. The sections $s,\,s'$ generate a
subline bundle $D$ of $E_{\mid C'}$ as above and the two sections
of $D$ restrict injectively to $z$. \cqfd {\it Last step.}  Lemma
\ref{smooth} shows that we have an isomorphism
$$ H^0(\widetilde W,r^*\mathcal{O}_{{\mathbb P}(H^0(E))_{curv}}(k+1))
\cong H^0(W'_U,pr_2^*\mathcal{O}_U(k+1)),$$
so that proposition \ref{tresdur} reduces to
\begin{prop} \label{finale}The map
$$pr_2^*:H^0(U,\mathcal{O}_U(k+1))\cong S^{k+1}H^0(S,E)^*\rightarrow
H^0(W'_U,pr_2^*(\mathcal{O}_U(k+1)))$$ is an isomorphism.
\end{prop}
{\bf Proof.} The proof works as in \cite{Vo} ; we note that $W'_U$
is the zero locus of a section $\tilde\sigma$ of a certain rank
$4$ vector bundle $pr_1^*\tilde E_2\otimes
pr_2^*(\mathcal{O}_U(1))$ on $\widetilde{S\times S}\times U$. We
use then the corresponding  Koszul resolution of
$\mathcal{I}_{W'_U}$ to conclude that
\begin{eqnarray}
\label{horrible} H^1(\widetilde{S\times S}\times U,{\mathcal
I}_{W'_U}\otimes pr_2^*(\mathcal{O}_U(k+1)))=0.
\end{eqnarray}
 There is one difference with the case considered in \cite{Vo}: namely, the spectral
 sequence which converges
 to
 $H^*(\widetilde{S\times S}\times U,{\mathcal I}_{W'_U}\otimes pr_2^*\mathcal{O}_U(k+1))$,
 has degree $1$ terms
 $$E_1^{i,1-i}=H^i(\widetilde{S\times S}\times U,pr_1^*\bigwedge^{i}\tilde E_2^*\otimes
  pr_2^*(\mathcal{O}_U(k+1-i))),\,i\geq1.$$
 Of course we have
 $$H^0(U,\mathcal{O}_U(k+1))=S^{k+1}H^0(S,E)^*.$$
 But unlike the case considered in \cite{Vo}, where we worked over the whole
 $\mathbb{P}(H^0(S,E))$, there might be
 some terms
 $$H^{i-1}(\widetilde{S\times S},\bigwedge^{i}\tilde E_2^*)\otimes H^1(U,\mathcal{O}_U(k+1))$$
 contributing to the
 term $E_1^{i,1-i}$ above. It turns out that this is not the case, thanks to the following
 lemma :
 \begin{lemm}We have
 $$H^{i-1}(\widetilde{S\times S},\bigwedge^{i}\tilde E_2^*)=0,$$ for any $i\geq1$.
 \end{lemm}
 {\bf Proof.}
 We refer to  \cite{Vo},  Proposition 6,
 for more details and similar computations.

 First of all, the vanishing $h^0(\widetilde{S\times S},\tilde
 E_2^*)=0$ follows from the fact that the dual vector bundle
 $\tilde E_2$ admits for space of global sections the space
 $H^0(E)$, which generates it generically, and that all of these
 sections vanish somewhere.

 Next, $\bigwedge^{4}\tilde E_2^*\cong det\,\tilde
 E_2^*=(-L)\boxtimes(-L)(2D)$ where $D$ is the
 exceptional divisor of  $\widetilde{S\times S}$, and we denote abusively  by
 $(-L)\boxtimes(-L)$ the pull-back of the line bundle $(-L)\boxtimes(-L)$
 on $S\times S$ to $\widetilde{S\times
 S}$
 via the blowing-down map.
  Since the
 canonical divisor of $\widetilde{S\times S}$ is equal to $D$,
 we get by Serre duality :
 $$ H^3(\widetilde{S\times S},\bigwedge^{4}\tilde E_2^*)=
 H^1(\widetilde{S\times S},L\boxtimes L(-D))^*,$$
 and the space on the right is $0$ because
 the multiplication map $S^2H^0(S,L)\rightarrow H^0(S,2L)$ is surjective.

  To compute $H^2(\widetilde{S\times S},\bigwedge^{3}\tilde E_2^*)$, we use the
  isomorphism
  $$\bigwedge^{3}\tilde E_2^*\cong det\,\tilde E_2^*\otimes
  \tilde E_2,$$
  and the exact sequence
  \begin{eqnarray}\label{nouvellexacte}
  0\rightarrow\tilde E_2\rightarrow pr_1^*E\oplus pr_2^*
  E\rightarrow \tau^*E\rightarrow0,
  \end{eqnarray}
  where $\tau:D\rightarrow S$ is the restriction of the
  blowing-down
  map to $D$, and the $pr_i$ are the projections to $S$ composed
  with the blowing-down map.
  It follows then from the associated long exact
  sequence that we only have  to prove the vanishings
  \begin{eqnarray}
  \label{ennuyeux1} H^1(D,\tau^*E(-2L)) (2D_{\mid D}))=0,\\
\label{ennuyeux2}
  H^2(\widetilde{S\times S}, pr_i^*E((-L)\boxtimes(-L))(2D))=0.
  \end{eqnarray}
  (\ref{ennuyeux1})  comes from the fact that
  $R^1\tau_* (2D_{\mid D})\cong\mathcal{O}_S$ and from
  $H^0(S,E(-2L))=0$.
  (\ref{ennuyeux2}) is deduced from Serre's duality and from
  $$H^2(\widetilde{S\times S},
  pr_i^*E^*(L\boxtimes L)(-D))=0.$$
  This last property is itself deduced from
  $H^2(S\times S,pr_i^*E^*(L\boxtimes L)=0$ and from
  $$H^1(D,pr_i^*E^*(L\boxtimes L)_{\mid
  D})=H^1(D,\tau^*(E(2L)))=0.$$

It remains only to prove the vanishing $H^1(\widetilde{S\times
S},\bigwedge^{2}\tilde E_2^* )=0$. We use for this the fact, which
follows from the dualization of the exact sequence
(\ref{nouvellexacte}) cf \cite{Vo}, that $\bigwedge^{2}\tilde
E_2^*$ has a filtration whose successive quotients are
$$\bigwedge^2(pr_1^*E^*\oplus pr_2^* E^*),\,
(pr_1^*E^*\oplus pr_2^*E^*)\otimes{\tau}^*E^*\otimes{\mathcal
O}_D(D),\, \bigwedge^2{\tau}^*E^*\otimes{\mathcal O}_D(2D).$$ It
is immediate to prove that each term has $H^1=0$.

 \cqfd

 Once we have these vanishings, the spectral sequence converging to
 $H^*(\widetilde{S\times S}\times U,{\mathcal I}_{W'_U}\otimes pr_2^*\mathcal{O}_U(k+1))$
 has the same shape in degree $1$ as in \cite{Vo},
 and then the proof  of the vanishing (\ref{horrible}) works as in \cite{Vo}.
 This concludes the proof of proposition
 \ref{finale}, hence of theorem \ref{evennode}.

\section{Proof of  Theorem \ref{ourthm}}
We start recalling the duality theorem of \cite{Gre}, which we
state here only in the case of surfaces :
\begin{theo}\label{dualthm}(Green) Let $X$ be a smooth projective surface,
$M$ be a line bundle on $X$ which is generated by sections, and
$\mathcal{F}$ be a coherent sheaf on $X$ satisfying the condition
$$H^1(X,\mathcal{F} (sM))=0,\,\forall
s\in\mathbb{Z}.$$ Then there is for all $p,\,q$ a duality
isomorphism (which is canonical up to a multiplicative
coefficient):
$$K_{p,q}(X,M,\mathcal{F})\cong K_{r-2-p,3-q}(X,M,\mathcal{F}^*
\otimes K_X)^*,$$ where $r+1=h^0(X,M)$.
\end{theo}

We consider now the case where $X$ is the $K3$ surface $S$ of the
previous section,  $M$ is either $L$ or $L'$, and $\mathcal{F}$ is
trivial. Then in the first case, $r+1=g+1=2k+2$, and in the second
case $r'+1=g'+1=2k+3$. So the duality theorem above gives, using
the fact that $K_S\cong \mathcal{O}_S$ :
\begin{eqnarray}\label{dual} K_{k,1}(S,L)^*\cong K_{k-1,2}(S,L),\\
\label{dual'} K_{k+1,1}(S,L')^*\cong K_{k-1,2}(S,L').
\end{eqnarray}
Theorem \ref{evennode} now says that $K_{k+1,1}(S,L')=0$ or
equivalently by (\ref{dual'})
\begin{eqnarray}
\label{vank+1} K_{k-1,2}(S,L')=0.
\end{eqnarray}
Next, recall that we want to prove that $K_{k,1}(S,L)=0$, and by
(\ref{dual}), this is equivalent to
\begin{eqnarray}\label{wanted}  K_{k-1,2}(S,L)=0.
\end{eqnarray}

Recalling that $L'=L+\Delta$, and choosing a $\sigma\in H^0(S,L')$
such that $\sigma$ generates $H^0(S,L')/H^0(S,L)$, (equivalently,
$\sigma$ nowhere vanishes along $\Delta$), we have now the
following :
\begin{prop}\label{generators} The space $K_{k-1,2}(S,L)$ is
generated as follows : consider the Koszul differential
\begin{eqnarray}\label{diff'}\delta':\bigwedge^{k-1}H^0(S,L)\otimes H^0(S,L-\Delta)\rightarrow
\bigwedge^{k-2}H^0(S,L)\otimes H^0(S,2L-\Delta).
\end{eqnarray}
For any $\alpha\in Ker\,\delta'$, multiplication on the right by
$\sigma\in H^0(S,L+\Delta)$ provides an element
$\alpha\cdot\sigma$ which is in
$$Ker\,\delta:\bigwedge^{k-1}H^0(S,L)\otimes H^0(S,2L)\rightarrow
\bigwedge^{k-2}H^0(S,L)\otimes H^0(S,3L),$$ where $\delta$ is also
the Koszul differential, however acting on a different space. The
classes of these elements $\alpha\cdot\sigma$ generate
$K_{k-1,2}(S,L)$.
\end{prop}
{\bf Proof.}  Let $\beta\in Ker\,\delta\subset
\bigwedge^{k-1}H^0(S,L)\otimes H^0(S,2L)$. Since we know that
$K_{k-1,2}(S,L')=0$ by (\ref{vank+1}), we can write
$$\beta=\delta\gamma$$
for some $\gamma\in \bigwedge^{k}H^0(S,L')\otimes H^0(S,L')$.
Since $H^0(S,L')=H^0(S,L)\oplus<\sigma>$, we can now decompose
$\gamma$ as
$$\gamma=\gamma_1+\sigma\wedge\gamma_2+\gamma_3\otimes\sigma+\sigma\wedge\gamma_4\otimes\sigma,$$
where
$$
\gamma_1\in\bigwedge^{k}H^0(S,L)\otimes H^0(S,L),\,\,\,\,
\gamma_2\in\bigwedge^{k-1}H^0(S,L)\otimes H^0(S,L),$$
$$\gamma_3\in\bigwedge^{k}H^0(S,L),\,\,\,\,
\gamma_4\in\bigwedge^{k-1}H^0(S,L) .$$ The fact that
$\delta\gamma=\beta$ belongs to $\bigwedge^{k-1}H^0(S,L)\otimes
H^0(S,2L)$ implies  that $\gamma_4=0$ since $\gamma_4$ identifies
to the image of $\delta\gamma$ in $\bigwedge^{k-1}H^0(S,L')\otimes
H^0(2L'_{\mid\Delta})$. Next, since we consider $\beta$ only
modulo
$$Im\,\delta:\bigwedge^{k}H^0(S,L)\otimes
H^0(S,L)\rightarrow\bigwedge^{k-1}H^0(S,L)\otimes H^0(S,2L),$$
 we may  assume, modifying $\beta$ by an exact element, that $\gamma_1=0$.
Finally, we note that $\gamma$ is defined up to $\delta$-closed
and in particular up to $\delta$-exact elements. Using the
relation
$$\gamma_3\otimes\sigma=-\delta(\sigma\wedge\gamma_3)-\sigma\wedge\delta\gamma_3,$$
we conclude that modifying $\gamma$ we may also assume that
$\gamma_3=0$.

In conclusion, $K_{k-1,2}(S,L)$ is generated by classes of
$\delta$-closed elements $\beta$ such that in
$\bigwedge^{k-1}H^0(S,L')\otimes H^0(S,2L')$, we have
$$\beta=\delta(\sigma\wedge\gamma),$$
for some $\gamma\in\bigwedge^{k-1}H^0(S,L)\otimes H^0(S,L).$ Now
we observe that the condition
$$\delta(\sigma\wedge\gamma)\in\bigwedge^{k-1}H^0(S,L)\otimes
H^0(S,2L)$$ implies that $\delta\gamma=0$. Hence
$$\delta(\sigma\wedge\gamma)=-\gamma\cdot\sigma.$$
The condition that
$$\gamma\cdot\sigma\in\bigwedge^{k-1}H^0(S,L)\otimes H^0(S,2L),$$
with $\sigma\in H^0(S,L+\Delta),\,\sigma_{\mid\Delta}\not=0$,
implies now that $\gamma\in\bigwedge^{k-1}H^0(S,L)\otimes
H^0(S,L-\Delta)$. Hence $\gamma\in Ker\,\delta'$ and the
proposition is proved. \cqfd Our next task is to compute the
dimension of the space
$$K:=Ker\,(\delta':\bigwedge^{k-1}H^0(S,L)\otimes H^0(S,L-\Delta)\rightarrow
\bigwedge^{k-2}H^0(S,L)\otimes H^0(S,2L-\Delta)).$$ Notice that
the Koszul complex of $(S,L,L-\Delta)$ equipped with the Koszul
differential $\delta'$ has the following shape :
\begin{eqnarray}\label{complex}
0\rightarrow\bigwedge^{k-1}H^0(S,L)\otimes
H^0(S,L-\Delta)\stackrel{\delta'}{\rightarrow}
\bigwedge^{k-2}H^0(S,L)\otimes H^0(S,2L-\Delta)\rightarrow\\
 \ldots \nonumber \rightarrow \bigwedge^{k-i}H^0(S,L)\otimes
H^0(S,iL-\Delta)\rightarrow\ldots.
\end{eqnarray}
So $K$ is the first cohomology group of this complex, while the
next ones are the  $K_{k-i,i-1}(S,L,L-\Delta)$ for $i\geq 2$. We
have now
\begin{lemm} The Koszul cohomology groups
$K_{k-i,i-1}(S,L,L-\Delta)$ vanish for $i\geq 2$.
\end{lemm}
{\bf Proof.} We observe that the triple $(S,L,L-\Delta)$ satisfies
the assumptions of the duality theorem \ref{dualthm}. Hence, using
$K_S\cong \mathcal{O}_S$, and $h^0(L)=2k+2$, we conclude that
$K_{k-i,i-1}(S,L,L-\Delta)$ is dual to
$K_{k-1+i,3-i+1}(S,L,-L+\Delta)$.

If $i=2$, the last group is the cohomology at the middle of the
sequence
$$ \bigwedge^{k+2}H^0(S,L)\otimes
H^0(S,\Delta)\stackrel{\delta_1}{\rightarrow}
\bigwedge^{k+1}H^0(S,L)\otimes
H^0(S,L+\Delta)\stackrel{\delta_2}{\rightarrow}
\bigwedge^{k}H^0(S,L)\otimes H^0(S,2L+\Delta).$$ Now we use the
equality $L'=L+\Delta$ and Theorem \ref{evennode} to conclude that
if $\beta\in Ker\,\delta_2$, then we have
$$\beta=\delta\gamma,$$
for some $\gamma\in \bigwedge^{k+2}H^0(S,L')$. As in the previous
proof, we now write
$$\gamma=\gamma_1+\sigma\wedge\gamma_2,$$
with
$$\gamma_1\in\bigwedge^{k+2}H^0(S,L),\,\gamma_2\in\bigwedge^{k+1}H^0(S,L).$$
The fact that
$$\delta\gamma=\beta\in\bigwedge^{k+1}H^0(S,L)\otimes
H^0(S,L+\Delta)$$ implies immediately that $\gamma_2=0$. Hence in
fact, we have in $\bigwedge^{k+1}H^0(S,L')\otimes H^0(S,L')$ the
equality
$$\beta=\delta\gamma,\,\gamma\in\bigwedge^{k+2}H^0(S,L).$$
Using the fact that the inclusion $H^0(S,L)\subset
H^0(S,L+\Delta)$ is the multiplication by the unique section of
$H^0(S,\Delta)$, it is obvious that this is equivalent to
$\beta\in Im\,\delta_1$. So the claim is proved in this case.

Next assume that $i=3$. Then $K_{k-1+i,3-i+1}(S,L,-L+\Delta)$ is
the cohomology  in the middle of the sequence
$$\bigwedge^{k+3}H^0(S,L)\otimes
H^0(S,-L+\Delta)\stackrel{\delta_1}{\rightarrow}
\bigwedge^{k+2}H^0(S,L)\otimes
H^0(S,\Delta)\stackrel{\delta_2}{\rightarrow}
\bigwedge^{k+1}H^0(S,L)\otimes H^0(S,L+\Delta).$$ But since
$H^0(S,\Delta)$ is of dimension $1$, it is easy to see that
$Ker\,\delta_2=0$. So this case is also proved.

Finally, if $i\geq 4$, $K_{k-1+i,3-i+1}(S,L,-L+\Delta)$ is $0$
because it is the cohomology at the middle of a complex with
vanishing middle term, since $H^0(S, sL+\Delta)=\{0\}$ for $s<0$.

\cqfd
\begin{coro}\label{dimension} The dimension of $K$ is equal to the binomial
coefficient $\begin{pmatrix}{2k+1}\\{k-1}\end{pmatrix}$.
\end{coro}
{\bf Proof.} $K$ is the degree $0$ cohomology group of the complex
(\ref{complex}) whose all next cohomology groups vanish. Hence the
dimension of $K$ is equal to the Euler characteristic of this
complex. Since the terms of the complex are
$\bigwedge^{k-i}H^0(S,L)\otimes H^0(S,iL-\Delta)$ put in degree
$i-1$, for $i\geq1$, and since
$$h^0(S,L)=2k+2,\,h^0(S,iL-\Delta)=1+2ki^2-2i,$$
we are reduced to proving the following identity:
$$\begin{pmatrix}{2k+1}\\{k-1}
\end{pmatrix}=\sum_{i\geq1}(-1)^{i-1}\begin{pmatrix}{2k+2}\\{k-i}\end{pmatrix}(1+2ki^2-2i).$$
The proof is left to the reader. \cqfd Recall now the vector
bundle $E$ from  (\ref{Lmuk}).  Our strategy to conclude the proof
of Theorem \ref{ourthm}, or equivalently the vanishing
$K_{k-1,2}(S,L)=0$, will be to construct a map
$$\phi:S^{k-1}H^0(S,E)\rightarrow K=Ker\,\delta'$$
and to prove first of all  that it is an isomorphism and secondly
that all the elements in $Im\,\phi$ are annihilated by the map
$\cdot\sigma$ of Proposition \ref{generators}. The vanishing
$K_{k-1,2}(S,L)=0$ will then be a consequence of Proposition
\ref{generators}.

\vspace{0,5cm}

 {\it Construction of $\phi$.} Recall  that
$E_{\mid\Delta}\cong\mathcal{O}_\Delta^2$, and that the
restriction map
$$H^0(S,E)\rightarrow H^0(\Delta,E_{\mid\Delta})$$
is surjective. Since $H^0(S,E)$ is of dimension $k+3$,
$H^0(S,E(-\Delta))$ is of dimension $k+1$. Consider the
determinant map
$$d:\bigwedge^2H^0(S,E)\rightarrow H^0(S,L').$$
 Note that for
$v\in\bigwedge^2H^0(S,E(-\Delta))$ we have $d(v)\in
H^0(S,L(-\Delta))$ and for $v\in H^0(S,E)\wedge
H^0(S,E(-\Delta))$, we have $d(v)\in H^0(S,L)$. Let
$$w_1,\ldots,w_{k+1}$$
be a basis of $H^0(S,E(-\Delta))$. The map $\phi$ is defined by
the following formula
\begin{eqnarray}\label{phi}\phi(\tau^{k-1})=\sum_{i<j}(-1)^{i+j}d(\tau\wedge
w_1)\wedge\ldots\hat i\ldots\hat j\ldots\wedge d(\tau\wedge
w_{k+1})\otimes d(w_i\wedge w_j).
\end{eqnarray}
By the remarks above, we have
$$\phi(\tau^{k-1})\in \bigwedge^{k-1}H^0(S,L)\otimes
H^0(S,L(-\Delta))\subset \bigwedge^{k-1}H^0(S,L')\otimes
H^0(S,L').$$ We prove now:
\begin{lemm} The image of $\phi$ is contained in $Ker\,\delta'$, where
$\delta'$ is the Koszul differential of (\ref{diff'}).
\end{lemm}
{\bf Proof.} Observe that we have the following quadratic
equations for $S$, imbedded in projective space via $\mid L'\mid$,
(these equations are in fact quadratic equations defining the
Grassmannian of codimension $2$ subspaces of $H^0(E)$, in which
$S$ lies naturally): consider the natural map
$$\psi:\bigwedge^3H^0(S,E)\otimes H^0(S,E)
\rightarrow \bigwedge^2H^0(S,E)\otimes \bigwedge^2H^0(S,E)$$
$$
\stackrel{d\otimes d}{\rightarrow}H^0(S,L')\otimes
H^0(S,L')\rightarrow S^2H^0(S,L').$$ Here the first map sends
$v_1\wedge v_2\wedge v_3\otimes\gamma$ to $$v_2\wedge v_3\otimes
v_1\wedge\gamma -v_1\wedge v_3\otimes v_2\wedge\gamma+v_1\wedge
v_2\otimes v_3\wedge\gamma.$$ We claim that the image of $\psi$ is
contained in the ideal of $S$. The reason is simply that the map
$\psi$ commutes with evaluation at $x\in S$ and that since
$rank\,E=2$, we have $\bigwedge^3E_x=0$.

So we conclude that we have the following equalities :
\begin{eqnarray}\label{equations}
d(v_2\wedge v_3)\cdot d(v_1\wedge\gamma) -d(v_1\wedge v_3)\cdot
d(v_2\wedge\gamma)\\ \nonumber +d(v_1\wedge v_2)\cdot d(
v_3\wedge\gamma)=0 \,\,\,\,{\rm in}\,\,\,\,H^0(S,2L').
\end{eqnarray}
We now compute :
$$\delta'(\phi(\tau^{k-1}))=\sum_{k<i<j}(-1)^{i+j+k}d(\tau\wedge
w_1)\wedge\ldots\hat k\,\hat i\,\hat j\,\ldots\wedge d(\tau\wedge
w_{k+1})\otimes d(\tau\wedge w_k)\cdot d(w_i\wedge w_j)$$
$$-\sum_{i<k<j}(-1)^{i+j+k}d(\tau\wedge
w_1)\wedge\ldots\hat i\,\hat k\,\hat j\ldots\wedge d(\tau\wedge
w_{k+1})\otimes d(\tau\wedge w_k)\cdot d(w_i\wedge w_j)$$
$$+\sum_{i<j<k}(-1)^{i+j+k}d(\tau\wedge
w_1)\wedge\ldots\hat i\,\hat j\,\hat k\ldots\wedge d(\tau\wedge
w_{k+1})\otimes d(\tau\wedge w_k)\cdot d(w_i\wedge w_j).$$ This is
also equal to
$$\sum_{i<j<k}(-1)^{i+j+k}d(\tau\wedge
w_1)\wedge\ldots\hat i\ldots\hat j\ldots\hat k\ldots\wedge
d(\tau\wedge w_{k+1})$$
$$\otimes (d(\tau\wedge w_k)\cdot d(w_i\wedge w_j)
+d(\tau\wedge w_i)\cdot d(w_j\wedge w_k)-d(\tau\wedge w_j)\cdot
d(w_i\wedge w_k)).$$ Hence by (\ref{equations}), we find that
$$\delta'(\phi(\tau^{k-1}))=0\,\,{\rm in}\,\,
\bigwedge^{k-2}H^0(S,L)\otimes H^0(S,2L-\Delta).$$ \cqfd
\begin{rema} {\rm The map $\phi$ is strongly related to the
construction due to Green and Lazarsfeld (see \cite{Gre},
Appendix) of non trivial syzygies in $K_{r_1+r_2-1,1}(X,L_1\otimes
L_2)$, where for $i=1,\,2$,  $L_i$ are line bundles on $X$ with
$r_i+1=h^0(X,L_i)$. The precise relation is obtained by taking
$X=C\in \mid L\mid$, $L_1$ a line bundle of degree $k+2$ on $C$
with $h^0(L_1)=2,\,h^0(L_1-\Delta_{\mid C})=1$, and $L_2=K_C-L_1$.
One has to use for that  the  relation (given by sequences like
(\ref{exE})) between the Lazarsfeld vector bundle $E$ and linear
systems on the curve $C$, or more precisely $C\cup\Delta$.}
\end{rema}
We shall prove  the following :
\begin{prop} \label{iso}The map
$$\phi:S^{k-1}H^0(S,E)\rightarrow K$$
is an isomorphism.
\end{prop}
{\bf Proof.} By corollary \ref{dimension}, both spaces have the
same dimension, since $rank\,H^0(S,E)=k+3$. The fact that $\phi $
is an isomorphism  reduces then  to the following :
\begin{prop}\label{injective}
The map $\phi$ is injective.
\end{prop}
We postpone the proof of Proposition \ref{injective} to the next
section. \cqfd

Assuming Proposition \ref{injective},  the proof of the vanishing
$$K_{k-1,2}(S,L)=0$$
is then a consequence of Proposition \ref{iso}, Proposition
\ref{generators} and of the following :
\begin{prop} For $\beta\in Im\,\phi\subset K$, we have
$$\beta\cdot\sigma=0\,\,\,{\rm in}\,\,\, K_{k-1,2}(S,L).$$
\end{prop}
{\bf Proof.}  Let $\beta=\phi(\tau^{k-1})$.  We may assume first
that $\tau\not\in H^0(S,E(-\Delta))$, and then  that
$\sigma=d(\tau\wedge w)$, for some $w\in H^0(S,E)$, because the
result depends only on the class of $\sigma$ modulo $H^0(S,L)$,
and the map
$$H^0(S,E)\rightarrow H^0(S,L+\Delta)/H^0(S,L),$$
$$v\mapsto d(\tau\wedge v)\,\,mod\,\,H^0(S,L)$$
is surjective. Next recall the formula (\ref{phi})
 $$\phi(\tau^{k-1})=\sum_{i<j}(-1)^{i+j}d(\tau\wedge
w_1)\wedge\ldots\hat i\ldots\hat j\ldots\wedge d(\tau\wedge
w_{k+1})\otimes d(w_i\wedge w_j).$$
 Using  equations (\ref{equations})
applied to $v_1=w_i,\,v_2=w_j,\,v_3=\tau,\,\gamma=w$, we get now
$$d(w_i\wedge w_j)\cdot d(\tau\wedge w)=d(\tau\wedge w_j)\cdot
d(w_i\wedge w) -d(\tau\wedge w_i)\cdot d(w_j\wedge w).$$ Hence
$$\phi(\tau^{k-1})\cdot d(\tau\wedge w)=\sum_{i<j}(-1)^{i+j}d(\tau\wedge
w_1)\wedge\ldots\hat i\ldots\hat j\ldots\wedge d(\tau\wedge
w_{k+1})$$
$$\otimes (-d(\tau\wedge w_i)\cdot d(w_j\wedge
w)+d(\tau\wedge w_j)\cdot d(w_i\wedge w)).$$ Now the  expression
on the right is equal to $\delta \beta'$, with
$$\beta'=\sum_i(-1)^i d(\tau\wedge w_1)\wedge\ldots\hat
i\ldots\wedge d(\tau\wedge w_{k+1})\otimes d(w_i\wedge w),$$ and
since $w_i\in H^0(S,E(-\Delta))$ we have
$$d(\tau\wedge w_i)\in H^0(S,L),\, d(w_i\wedge w)\in H^0(S,L)$$
so that $\beta'\in \bigwedge^kH^0(S,L)\otimes H^0(S,L)$. So
$\beta\cdot\sigma=0$ in $K_{k-1,2}(S,L)$. \cqfd

\section{\label{sec3}Proof of Proposition \ref{injective}}
Let us first recall  the statement: we have the determinant map
$$d:\bigwedge^2H^0(E)\rightarrow H^0(S,L'),$$
which has the property that it does not vanish on non-zero
decomposable elements. Here the rank of $H^0(E)$ is $k+3$ and the
rank of $H^0(E(-\Delta))$ is $k+1$. We defined the map
$$\phi:S^{k-1}H^0(E)\rightarrow \bigwedge^{k-1}H^0(L)\otimes
H^0(L(-\Delta)),$$  explicitly by the formula (cf (\ref{phi}))
$$\phi(\tau^{k-1})=\sum_{i<j}(-1)^{i+j}d(\tau\wedge
w_1)\wedge\ldots\hat i\ldots\hat j\ldots\wedge d(\tau\wedge
w_{k+1})\otimes d(w_i\wedge w_j),$$ where the $w_l$'s form a basis
of $H^0(E(-\Delta))$. Proposition \ref{injective} states that this
map is injective.

 We give an ad hoc,
presumably not optimal, proof of this, relying on the particular
geometry of the determinant map $d$. We believe that it is in fact
true for any $d$ satisfying the condition that $d$ does not vanish
on decomposable elements.

We assume in the following that $k\geq2$.

In our situation, let $x\in S$ be a generic point. Consider the
composition $\phi_x$ of $\phi$ with the evaluation at $x$:
$$\phi_x:S^{k-1}H^0(E)\rightarrow\bigwedge^{k-1}H^0(S,L)\otimes
H^0(S,L(-\Delta)_{\mid x}).$$ Choose the basis $w_1,\ldots
w_{k+1}$ in such a way that $w_1,\dots,w_{k-1}$ form a basis of
$H^0(S,E(-\Delta)\otimes {\mathcal   I}_x)$. Then the $d(w_i\wedge
w_j)$ vanish at $x$ if $i$ or $j$ is non greater than $k-1$, while
$d(w_k\wedge w_{k+1})$ does not vanish in $H^0(S,L(-\Delta)_{\mid
x})$. Identifying this last space with $\mathbb{C}$, it follows
that $\phi_x$ has the following form up to a  coefficient:
\begin{eqnarray}
\label{mardi29}\phi_x(\tau^{k-1})=d(\tau\wedge
w_1)\wedge\ldots\wedge d(\tau\wedge w_{k-1}).
\end{eqnarray}

\vspace{0,5cm}

{\it First step.} We first use formula (\ref{mardi29}) to express
the map $\phi$, or rather its transpose, as the map induced in
cohomology by the top exterior power of a vector bundle map over
an adequate variety. That will allow us later on to use the Koszul
resolution of such top exterior powers.

Denote by ${\mathcal V}$ the vector bundle on the $K3$ surface,
which is defined by the exact sequence
$$0\rightarrow {\mathcal
V}\rightarrow H^0(E(-\Delta))\otimes{\mathcal   O}_S\rightarrow
E(-\Delta)\rightarrow0.$$ So the fiber $\mathcal{V}_x$ at $x\in S$
is the space generated by the $w_1,\ldots,w_{k-1}$ introduced
above.

 On $Y:=\mathbb{P}(H^0(E))\times S$, there is a natural map
$$h:pr_2^*{\mathcal   V}\otimes
pr_1^*{\mathcal   O}_{\mathbb{P}(H^0(E))}(-1) \rightarrow
H^0(S,L)\otimes{\mathcal   O}_Y,$$ which at the point $(\tau,x)$ is
the map
$$d(\tau\wedge\,\,\,):H^0(E(-\Delta)\otimes{\mathcal   I}_x)\rightarrow
H^0(L).$$ This map is injective when $\tau\not\in
H^0(E(-\Delta)\otimes{\mathcal   I}_x)$, and has for kernel
$<\tau>$ otherwise. Since we want to study the map induced in
cohomology by the top exterior power of $h$, we first want to make
$h$ into a morphism which is everywhere injective. This is done as
follows  (we refer to diagram (\ref{diagram1}) for the notations)
: Let
$$\mathbb{P}(\mathcal{V})=:Z=\{(\tau,x)\in Y,\,\tau\in
H^0(E(-\Delta)\otimes{\mathcal   I}_x)\}.$$ So $Z$ is the locus of
points $(\tau,x)$ where $h_{\tau,x}$ is not injective. Denote by
$f:\widetilde{Y}\rightarrow Y$ the blow-up of $Y$ along $Z$. For
simplicity, denote by $\mathbb{P}$ the space $\mathbb{P}(H^0(E))$
and  by $p$ the map $pr_1\circ f:\widetilde{Y}\rightarrow
\mathbb{P}$. Let also $q:=pr_2\circ f:\widetilde{Y}\rightarrow S$.
 The map $h$ extends
to a map
$$\tilde h:{\mathcal   G}\rightarrow H^0(S,L)\otimes
{\mathcal   O}_{\widetilde{Y}},$$ which is now injective everywhere,
where ${\mathcal   G}$ is obtained from $q^*{\mathcal V}\otimes
p^*{\mathcal   O}_{\mathbb{P}}(-1)$ by an elementary transform along
 the exceptional divisor $D$ of $f$. Namely ${\mathcal   G}$ fits in
an exact sequence
$$0\rightarrow q^*{\mathcal   V}\otimes
p^*{\mathcal   O}_{\mathbb{P}}(-1)\rightarrow {\mathcal
G}\rightarrow{\mathcal   H}_D\rightarrow0,$$ where ${\mathcal
H}_D$ is a line bundle supported on $D$, and the restriction of
the first map to $D$ has exactly for kernel the kernel of the map
$h_{\mid D}$, that is the sub-line bundle
$$p^*{\mathcal   O}_{\mathbb{P}}(-2)_{\mid D}\subset q^*{\mathcal   V}\otimes
p^*{\mathcal   O}_{\mathbb{P}}(-1)_{\mid D}.$$

 Note that
 \begin{eqnarray}
 \label{formuledet}\bigwedge^{k-1}{\mathcal   G}=det\,{\mathcal   G}=
 p^*{\mathcal   O}_{\mathbb{P}}(-k+1)\otimes
 q^*(L^{-1}(\Delta))(D).
 \end{eqnarray}

 Let
 $$\tilde{h}_{k-1}:det\,{\mathcal   G}\rightarrow\bigwedge^{k-1}H^0(L)\otimes{\mathcal   O}_{\widetilde{Y}}
 $$ be the map which is the $k-1$-th exterior power of
 $\tilde{h}$, and let
 $$h'_{k-1}:\bigwedge^{k-1}H^0(L)^*\otimes
 q^*(L^{-1}(\Delta))(D)\rightarrow \bigwedge^{k-1}{\mathcal G}^*\otimes
 q^*(L^{-1}(\Delta))(D)$$
 be the transpose of $\tilde{h}_{k-1}$ twisted by
 $q^*(L^{-1}(\Delta))(D)$.
We first claim
  that the transpose of the map $\phi$ identifies
 to the map $h^{2}({h}'_{k-1})$:
  \begin{eqnarray}\label{h2k}
  \begin{matrix}&H^{2}( \widetilde{Y}, \bigwedge^{k-1}H^0(L)^*\otimes
 q^*(L^{-1}(\Delta))(D))&
  \rightarrow&
  H^{2}(\widetilde{Y},\bigwedge^{k-1}{\mathcal G}^*\otimes q^*(L^{-1}(\Delta))(D) )&\\
  &\parallel&&\parallel&\\
  &\bigwedge^{k-1}H^0(L)^*\otimes H^0(L(-\Delta))&\stackrel{^t\phi}{\rightarrow}&
  H^2(\widetilde{Y},p^*\mathcal{O}(k-1))=S^{k-1}H^0(E)^*&.
  \end{matrix}
\end{eqnarray}
  Indeed, by (\ref{formuledet}),
  we have
  $$R^{k+2}q_*(K_{\widetilde Y}\otimes p^*{\mathcal   O}_{\mathbb{P}}(-k+1))=
R^{k+2}q_*(K_{\widetilde Y}\otimes det\,{\mathcal   G}\otimes
  q^*(L(-\Delta))(-D))
  ,$$
and by Serre duality and $K_S=\mathcal{O}_S$, the left hand side
identifies to $S^{k-1}H^0(E)\otimes{\mathcal O}_S$.
  Now, formula (\ref{mardi29}) says that
  the map induced (up to a twist) by $ \tilde{h}_{k-1}$ :
  $$S^{k-1}H^0(E)\otimes{\mathcal O}_S\cong
  R^{k+2}q_*(K_{\widetilde Y}\otimes det\,{\mathcal   G}\otimes
  q^*(L(-\Delta))(-D))$$
  $$\rightarrow
R^{k+2}q_*(K_{\widetilde{Y}}\otimes\bigwedge^{k-1}H^0(L)\otimes
q^*(L(-\Delta))(-D))= \bigwedge^{k-1}H^0(L)\otimes L(-\Delta)$$ is
exactly the map $\phi$ followed with evaluation.
 Taking global sections on $S$, we conclude that
 $\phi$ is the map induced by $ \tilde{h}_{k-1}$ (up to a twist):
 $$H^{k+2}(\widetilde{Y},K_{\widetilde Y}\otimes det\,{\mathcal   G}\otimes
  q^*(L(-\Delta))(-D))\rightarrow H^{k+2}(\widetilde{Y},K_{\widetilde{Y}}\otimes\bigwedge^{k-1}H^0(L)\otimes
q^*(L(-\Delta))(-D)),$$
 and applying Serre duality gives
 the result.

  So the content of Proposition \ref{injective} is the surjectivity of
  the map $h^{2}({h}'_{k-1})$.

\vspace{0,5cm}

{\it Second step.} We shall now analyse the spectral sequence
associated to the Koszul resolution  of $Ker\,h'_{k-1}$:
   associated to the surjective map
  $^t\tilde{h}$, there is a resolution
  \begin{eqnarray}\label{com}
  0\rightarrow  \bigwedge^{2k+2}H^0(L)^*\otimes
  S^{k+3}{\mathcal   G}\rightarrow\ldots\rightarrow
\bigwedge^{k}H^0(L)^*\otimes
  {\mathcal   G}\\\rightarrow
 \bigwedge^{k-1}H^0(L)^*\otimes{\mathcal   O}_{\widetilde{Y}}
\rightarrow\bigwedge^{k-1}{\mathcal   G}^*\rightarrow0. \nonumber
  \end{eqnarray}
   We claim now that the surjectivity of the map (\ref{h2k})
   follows  from the following
  lemmas.
  \begin{lemm}\label{lem5mai}For $1\leq l\leq k$, or $l=k+2$, we have
$$H^{l+2}(\widetilde{Y},S^l{\mathcal G}\otimes q^*(L^{-1}(\Delta))(D))=0.$$
  \end{lemm}
  \begin{lemm}\label{le2mardi29}
  The sequence
  \begin{eqnarray}\label{complexe9mai} \bigwedge^{2k+1}
  H^0(L)^*\otimes H^{k+3}(\widetilde{Y},S^{k+2}{\mathcal   G}\otimes
  q^*(L^{-1}(\Delta))(D))\\
  \nonumber
\rightarrow
  \bigwedge^{2k}H^0(L)^*\otimes H^{k+3}(\widetilde{Y},S^{k+1}{\mathcal   G}\otimes
  q^*(L^{-1}(\Delta))(D))\\
  \nonumber
\rightarrow
  \bigwedge^{2k-1}H^0(L)^*\otimes H^{k+3}(\widetilde{Y},S^{k}{\mathcal   G}\otimes
  q^*(L^{-1}(\Delta))(D))
  \end{eqnarray}
  induced by the complex (\ref{com}) is
  exact at the middle.
  \end{lemm}

Indeed,  these two lemmas together imply that the vector bundle
$Ker\,{h}'_{k-1}$ satisfies $H^3(\widetilde{Y},Ker\,h'_{k-1})=0$,
which implies the surjectivity of $h^2(h'_{k-1})$: In fact,
twisting the complex (\ref{com}) by
  $q^*(L^{-1}(\Delta))(D)$, we get a resolution of $Ker\,{h}'_{k-1}$ as follows :
$$0\rightarrow  \bigwedge^{2k+2}H^0(L)^*\otimes
  S^{k+3}{\mathcal   G}\otimes q^*(L^{-1}(\Delta))(D)\rightarrow\ldots\rightarrow
\bigwedge^{k}H^0(L)^*\otimes
  {\mathcal   G}\otimes q^*(L^{-1}(\Delta))(D) $$
  $$\rightarrow
  Ker\,{h}'_{k-1}
 \rightarrow0.
$$
The associated spectral sequence abutting to the hypercohomology
of the complex $$0\rightarrow  \bigwedge^{2k+2}H^0(L)^*\otimes
  S^{k+3}{\mathcal   G}\otimes q^*(L^{-1}(\Delta))(D)\rightarrow\ldots\rightarrow
\bigwedge^{k}H^0(L)^*\otimes
  {\mathcal   G}\otimes q^*(L^{-1}(\Delta))(D) \rightarrow0,$$ where
  we put the last  term on the right  in degree $0$,
   has
   $$E_1^{p,q}\cong \bigwedge^{k-p}H^0(L)^*\otimes
   H^q(\widetilde{Y}, S^{1-p} {\mathcal   G}\otimes
   q^*(L^{-1}(\Delta))(D)).$$
Now Lemma \ref{lem5mai} says that the $E_1^{p,q}$ for
$p+q=3,\,q\geq3$ are $0$ unless $q=k+3$. For $q=k+3$, we have
$E_1^{-k,k+3}\not=0$ but  Lemma \ref{le2mardi29} says that
$E_2^{-k,k+3}=0$. Hence this complex has $\mathbb{H}^3=0$ and thus
$H^3(\widetilde{Y},Ker\,h'_{k-1})=0$.
 \cqfd

 \vspace{0,5cm}

 {\it Third step.} We start now proving Lemmas \ref{lem5mai} and
 \ref{le2mardi29}. We shall use for this  another geometric definition of
 the vector bundle $\mathcal{G}$. We refer to diagram (\ref{diagram1}) for
 the notations.

\vspace{0,5cm}

  {\bf Proof of Lemma \ref{lem5mai}.}
  Let ${\mathcal F}$ be the quotient bundle
$H^0(S,E)\otimes {\mathcal O}_S/{\mathcal   V}$. There is the
relative projection
$$\chi:\widetilde{Y}\rightarrow \mathbb{P}({\mathcal F}),$$
which makes $\widetilde{Y}$ isomorphic to $\mathbb{P}({\mathcal
H})$, where ${\mathcal H}$ is a vector bundle on
$\mathbb{P}({\mathcal   F})$ which fits in the exact sequence
\begin{eqnarray}
\label{f5mai} O\rightarrow \pi^*{\mathcal   V}\rightarrow{\mathcal
H}\rightarrow {\mathcal   O}_{\mathbb{P}({\mathcal
F})}(-1)\rightarrow0, \end{eqnarray} where $\pi:\mathbb{P}({\mathcal
F})\rightarrow S$ is the structural map.
 We
observe now that ${\mathcal G}$ is naturally isomorphic to the
twisted relative tangent bundle $T_\chi\otimes p^*{\mathcal
O}_{\mathbb{P}}(-2)$. To see this, we consider the relative Euler
sequence \begin{eqnarray} \label{Eurel}0\rightarrow p^*{\mathcal
O}_{\mathbb{P}}(-2)\rightarrow \chi^*{\mathcal   H}\otimes
p^*{\mathcal O}_{\mathbb{P}}(-1)\rightarrow T_\chi(-2)\rightarrow0.
\end{eqnarray}
 It induces a map
$$q^*{\mathcal   V}\otimes
p^*{\mathcal   O}_{\mathbb{P}}(-1)\rightarrow T_\chi,$$ and it
follows from the Euler sequence  that this map is injective away
from $D$ and has $p^*{\mathcal   O}_{\mathbb{P}}(-2)_{\mid D}$ as
kernel along $D$. Hence $T_\chi$ is deduced from $q^*{\mathcal
V}\otimes p^*{\mathcal   O}_{\mathbb{P}}(-1)$ by the same elementary
transform as ${\mathcal G}$.

\begin{eqnarray}\begin{diagram}\label{diagram1}
\node{D}\arrow{e,J}\arrow{s}\node{\widetilde{Y}=\mathbb{P}(\mathcal{H})}
\arrow{ssw,l,3}{p}\arrow{s,r}{f}\arrow{sse,r}{q}\arrow[2]{e,t}{\chi}\node[2]{\mathbb{P}(\mathcal{F})}
\arrow{ssw,r}{\pi}
\\ \node{Z=\mathbb{P}(\mathcal{V})}\arrow{e,J}
\node{Y}\arrow{sw,r}{pr_1}\arrow{se,r}{pr_2}\node{}\node{}
\\ \node{\mathbb{P}:=\mathbb{P}(H^0(E))}\node{}\node{S}
\end{diagram}
\end{eqnarray} The
relative Euler sequence (\ref{Eurel}) describes ${\mathcal G}$ by
the exact sequence :
\begin{eqnarray}\label{3maiexact}
0\rightarrow p^*{\mathcal O}_{\mathbb P}(-2) \rightarrow
 p^*{\mathcal O}_{\mathbb P}(-1)\otimes \chi^*{\mathcal H}
\rightarrow {\mathcal G} \rightarrow0.
\end{eqnarray}
Taking the $l$-th symmetric power,  we get the exact sequence :
\begin{eqnarray}\label{3maiexactSl}
0 \rightarrow
 p^*{\mathcal O}_{\mathbb P}(-l-1)\otimes \chi^*S^{l-1}{\mathcal H}
\rightarrow
 p^*{\mathcal O}_{\mathbb P}(-l)\otimes \chi^*{S^l\mathcal H}
\rightarrow S^l{\mathcal G} \rightarrow0.
\end{eqnarray}
Assume first that $1\leq l\leq k-1$. By the exact sequence
(\ref{3maiexactSl}),
 the vanishing
$$H^{l+2}(\widetilde{Y}, S^l{\mathcal G}\otimes q^*(L^{-1}(\Delta))(D))=0$$
is implied by the vanishings :
\begin{eqnarray}\label{van3mai}
H^{l+2}(\widetilde{Y}, p^*{\mathcal O}_{\mathbb P}(-l)\otimes
\chi^*{S^l\mathcal H}\otimes q^*(L^{-1}(\Delta))(D))=0
,\\
\nonumber
 H^{l+3}(\widetilde{Y}, p^*{\mathcal O}_{\mathbb P}(-l-1)\otimes
\chi^*{S^{l-1}\mathcal H}\otimes
q^*(L^{-1}(\Delta))(D))=0.\end{eqnarray} Now, if $k-1\geq l\geq2$,
the line bundles $p^*{\mathcal O}_{\mathbb
P}(-l)(D),\,p^*{\mathcal O}_{\mathbb P}(-l-1)(D)$, have trivial
cohomology  along the fibers of $\chi$, which are ${\mathbb
P}^{k-1}$'s, on which ${ D}$ restricts to ${\mathcal O}(1)$. Hence
the vanishings (\ref{van3mai}) are proved in this case. The case
$l=1$ is also easy.

If $l=k$, the argument above gives an inclusion
$$H^{k+2}(\widetilde{Y}, S^{k}{\mathcal G}\otimes
q^*(L^{-1}(\Delta))(D))$$
$$\hookrightarrow
H^{k+3}(\widetilde{Y}, p^*{\mathcal O}_{\mathbb P}(-k-1)\otimes
\chi^*S^{k-1}{\mathcal H}\otimes q^*(L^{-1}(\Delta))(D)).$$ By Serre
duality, this dualizes as
$$H^{1}(\widetilde{Y}, p^*{\mathcal O}_{\mathbb P}(-2)\otimes
\chi^*S^{k-1}{\mathcal H}^*\otimes q^*(L(-\Delta))(2D)).$$ Since
$p^*{\mathcal O}_{\mathbb
P}(-2)(2D)=\chi^*\mathcal{O}_{\mathbb{P}(\mathcal{F})}(-2)$ and $
S^{k-1}{\mathcal H}^*=R^0\chi_*(p^*{\mathcal O}_{\mathbb
P}(k-1))$, the last space is equal to
$$H^{1}(\widetilde{Y},p^*{\mathcal O}_{\mathbb
P}(k-3)(2D)\otimes q^*(L-(\Delta))),$$ which is  $0$ because the
map $f=(p,q)$ is the blow-down map, and $2D$ has trivial
cohomology along the fibers of $f$. It follows that
$H^{1}(\widetilde{Y},p^*{\mathcal O}_{\mathbb P}(k-3)(2D)\otimes
q^*(L-(\Delta)))=H^1(Y,pr_1^*{\mathcal O}_{\mathbb P}(k-3)\otimes
pr_2^*(L-(\Delta)))=0$.

 Next, if $l=k+2$, we use  the inclusion
$p^*{\mathcal O}_{\mathbb P}(-1)\otimes q^*{\mathcal V} \subset
{\mathcal G} $, which is an isomorphism away from $D$, to get a
surjective map
$$H^{k+4}(\widetilde{Y},p^*{\mathcal O}_{\mathbb P}(-k-2)
\otimes q^*S^{k+2}{\mathcal V}\otimes q^*(L^{-1}(\Delta))(D))$$
$$\rightarrow H^{k+4}(\widetilde{Y}, S^{k+2}{\mathcal G}\otimes
q^*(L^{-1}(\Delta))(D)).$$ The left hand side is zero because $D$
has trivial cohomology along the fibers of $(p,q)$, so that
$$H^{k+4}(\widetilde{Y},p^*{\mathcal O}_{\mathbb P}(-k-2) \otimes
q^*S^{k+2}{\mathcal V}\otimes q^*(L^{-1}(\Delta))(D))$$ $$=
H^{k+4}(Y,pr_1^*{\mathcal O}_{\mathbb P}(-k-2)\otimes
pr_2^*(S^{k+2}{\mathcal V}\otimes L^{-1}(\Delta)))=0$$ (recall
that $Y=\mathbb{P}(H^0(E))\times S$, with $rk\,H^0(E)=k+3$). Hence
we conclude that $ H^{k+4}(\widetilde{Y}, S^{k+2}{\mathcal
G}\otimes q^*(L^{-1}(\Delta))(D))=0$.
  \cqfd

 In order to prove Lemma \ref{le2mardi29}, we will need the
  following :
  \begin{lemm}
  \label{le1mardi29}\begin{enumerate}
\item \label{item1}The space
\begin{eqnarray}\label{group8mai}H^{k+3}(\widetilde{Y},S^{k+1}{\mathcal G}\otimes
q^*(L^{-1}(\Delta))(D))
\end{eqnarray}
is canonically isomorphic to $S^{k-2}H^0(E)\otimes H^2(S,{\mathcal
V}\otimes L^{-1}(\Delta))$.

\item \label{item2} The space
$H^{k+3}(\widetilde{Y},S^{k+2}{\mathcal G}\otimes
q^*(L^{-1}(\Delta))(D))$ is canonically isomorphic to
$S^{k-1}H^0(E)\otimes H^2(S,S^2{\mathcal V}\otimes
L^{-1}(\Delta))$.

\item \label{item3} The space
$H^{k+3}(\widetilde{Y},S^{k}{\mathcal G}\otimes
q^*(L^{-1}(\Delta))(D))$ is canonically isomorphic to
$S^{k-3}H^0(E)\otimes H^2(S, L^{-1}(\Delta))$.

\end{enumerate}
\end{lemm}

{\bf Proof.}
 We use the exact sequence
$$0\rightarrow \chi^*S^{l-1}{\mathcal H}
\otimes p^*{\mathcal O}_{\mathbb P}(-l-1) \rightarrow
\chi^*S^{l}{\mathcal H} \otimes p^*{\mathcal O}_{\mathbb P}(-l)
\rightarrow S^{l}{\mathcal G} \rightarrow0.$$ It implies by the
associated long exact sequence that the space $$
H^{k+3}(\widetilde{Y},S^{l}{\mathcal G}\otimes
q^*(L^{-1}(\Delta))(D))$$is isomorphic to
$$Ker\,(H^{k+4}(\widetilde{Y},
\chi^*S^{l-1}{\mathcal H} \otimes p^*{\mathcal O}_{\mathbb
P}(-l-1)\otimes q^*(L^{-1}(\Delta))(D))$$
$$\rightarrow H^{k+4}(\widetilde{Y},
\chi^*S^{l}{\mathcal H} \otimes p^*{\mathcal O}_{\mathbb
P}(-l)\otimes q^*(L^{-1}(\Delta))(D))).$$ Recalling that $$
p^*\mathcal{O}_{\mathbb{P}}(1)=\chi^*\mathcal{O}_{\mathbb{P}(\mathcal{F})}(1)(D),$$
we can rewrite this as \begin{eqnarray} \label{le89juin1}
Ker\,(H^{k+4}(\widetilde{Y}, \chi^*S^{l-1}{\mathcal H}
\otimes\chi^*\mathcal{O}_{\mathbb{P}(\mathcal{F})}(-1)\otimes
p^*{\mathcal O}_{\mathbb P}(-l)\otimes q^*(L^{-1}(\Delta)))\\
\nonumber \rightarrow H^{k+4}(\widetilde{Y}, \chi^*S^{l}{\mathcal
H}\otimes \chi^*\mathcal{O}_{\mathbb{P}(\mathcal{F})}(-1) \otimes
p^*{\mathcal O}_{\mathbb P}(-l+1)\otimes q^*(L^{-1}(\Delta)))).
\end{eqnarray} Now, we use the formula
$$K_{\widetilde{Y}/\mathbb{P}(\mathcal{F})}=p^*{\mathcal O}_{\mathbb
P}(-k)\otimes \chi^*det\,\mathcal{H}^*,$$ that is, by the exact
sequence (\ref{f5mai}), which gives
$det\,\mathcal{H}^*=\pi^*(L(-\Delta))\otimes
\mathcal{O}_{\mathbb{P}(\mathcal{F})}(1)$,
$$K_{\widetilde{Y}/\mathbb{P}(\mathcal{F})}=p^*{\mathcal O}_{\mathbb
P}(-l)\otimes
\chi^*\mathcal{O}_{\mathbb{P}(\mathcal{F})(1)}\otimes
q^*(L(-\Delta)).$$ It follows then by Leray spectral sequence and
relative Serre duality, that  (\ref{le89juin1})  is also equal to:
$$Ker\,( H^5(\mathbb{P}(\mathcal{F}),S^{l-1}\mathcal{H}\otimes S^{l-k}\mathcal{H}\otimes
\mathcal{O}_{\mathbb{P}(\mathcal{F})}(-2)\otimes
\pi^*(L^{-2}(2\Delta)))
  $$
$$\rightarrow  H^5(\mathbb{P}(\mathcal{F}),S^{l}\mathcal{H}\otimes S^{l-k-1}\mathcal{H}\otimes
\mathcal{O}_{\mathbb{P}(\mathcal{F})}(-2)\otimes
\pi^*(L^{-2}(2\Delta))  )),$$ where we make the convention that
negative symmetric powers are $0$, and where the map is induced by
the natural map
$$S^{l-1}\mathcal{H}\otimes S^{l-k}\mathcal{H}
\rightarrow S^{l}\mathcal{H}\otimes S^{l-k-1}\mathcal{H}.$$ We now
apply again relative Serre duality and Leray spectral sequence to
conclude that this is also equal to :

$$Ker\,(H^{k+4}(\widetilde{Y},
  p^*{\mathcal O}_{\mathbb P}(-k-l)\otimes \chi^*S^{l-k}{\mathcal H}
\otimes q^*(L^{-1}(\Delta))(D))$$
$$\rightarrow H^{k+4}(\widetilde{Y},
p^*{\mathcal O}_{\mathbb P}(-k-l-1)\otimes
\chi^*S^{l-k-1}{\mathcal H} \otimes  q^*(L^{-1}(\Delta))(D))).$$

We now distinguish according to the value of $l$.

\vspace{0,5cm}

 - {\it Case $l=k$}. We proved that in this case we
have
$$H^{k+3}(\widetilde{Y},S^{k}{\mathcal G}\otimes
q^*(L^{-1}(\Delta))(D))\cong H^{k+4}(\widetilde{Y},
  p^*{\mathcal O}_{\mathbb P}(-2k)
\otimes q^*(L^{-1}(\Delta))(D)).$$ Since
$K_{\widetilde{Y}}=p^*\mathcal{O}_{\mathbb{P}}(-k-3)(3D)$, this is
Serre dual to
$$ H^0(\widetilde{Y}, p^*{\mathcal O}_{\mathbb P}(k-3)
\otimes q^*(L(-\Delta))(2D))=S^{k-3}H^0(E)^*\otimes
H^0(S,L(-\Delta)).$$ Applying Serre's duality on $S$ gives then
\ref{item3}.

\vspace{0,5cm}

- {\it Case $l=k+1$}. In this case, we have
\begin{eqnarray}
\label{point1}H^{k+3}(\widetilde{Y},S^{k+1}{\mathcal G}\otimes
q^*(L^{-1}(\Delta))(D))\cong\\ \nonumber
Ker\,(H^{k+4}(\widetilde{Y},
  p^*{\mathcal O}_{\mathbb P}(-2k-1)\otimes\chi^*\mathcal{H}
\otimes q^*(L^{-1}(\Delta))(D)) \\ \nonumber \rightarrow
H^{k+4}(\widetilde{Y},
  p^*{\mathcal O}_{\mathbb P}(-2k-2)
\otimes q^*(L^{-1}(\Delta))(D))). \end{eqnarray} The right hand
side is computed as before : by Serre duality on $\widetilde{Y}$,
we get {\begin{eqnarray}\label{point2}H^{k+4}(\widetilde{Y},
  p^*{\mathcal O}_{\mathbb P}(-2k-2)
\otimes q^*(L^{-1}(\Delta))(D))\cong S^{k-1}H^0(E)\otimes
H^2(S,L^{-1}(\Delta)).
\end{eqnarray}

To compute the left hand side, we use the exact sequence
(\ref{f5mai}): $$ O\rightarrow \pi^*{\mathcal
V}\rightarrow{\mathcal H}\rightarrow {\mathcal
O}_{\mathbb{P}({\mathcal F})}(-1)\rightarrow0.$$ Pulling back to
$\widetilde{Y}$ and tensoring with $p^*{\mathcal O}_{\mathbb
P}(-2k-1) \otimes q^*(L^{-1}(\Delta))(D)$, it provides the exact
sequence: \begin{eqnarray} \label{le8ex2}O\rightarrow q^*{\mathcal
V}\otimes p^*{\mathcal O}_{\mathbb P}(-2k-1) \otimes
q^*(L^{-1}(\Delta))(D)\\ \nonumber \rightarrow\chi^*{\mathcal
H}\otimes p^*{\mathcal O}_{\mathbb P}(-2k-1) \otimes
q^*(L^{-1}(\Delta))(D)\\ \nonumber \rightarrow \chi^*{\mathcal
O}_{\mathbb{P}({\mathcal F})}(-1)\otimes p^*{\mathcal O}_{\mathbb
P}(-2k-1) \otimes q^*(L^{-1}(\Delta))(D)\rightarrow0.
\end{eqnarray} Note that, as (\ref{f5mai})
is (non-canonically) split along the fibers of $\pi$, the sequence
(\ref{le8ex2}) is (non-canonically) split along the fibers of $q$.
Hence there is an induced exact sequence on $S$:
\begin{eqnarray} \label{le8ex3}O\rightarrow {\mathcal
V}\otimes L^{-1}(\Delta)\otimes R^{k+2}q_* ( p^*{\mathcal
O}_{\mathbb P}(-2k-1)  (D))\\ \nonumber\rightarrow
L^{-1}(\Delta)\otimes R^{k+2}q_* (\chi^*{\mathcal H}\otimes
p^*{\mathcal O}_{\mathbb P}(-2k-1) (D))\\ \nonumber \rightarrow
L^{-1}(\Delta)\otimes R^{k+2}q_* ( \chi^*{\mathcal
O}_{\mathbb{P}({\mathcal F})}(-1)\otimes p^*{\mathcal O}_{\mathbb
P}(-2k-1)  (D))\rightarrow0.
\end{eqnarray}
Using the equality
$$K_{\widetilde{Y}}=K_{\widetilde{Y}/S}=p^*\mathcal{O}_{\mathbb{P}}(-k-3)(3D)$$
and relative Serre duality, we get :
$$R^{k+2}q_* ( p^*{\mathcal
O}_{\mathbb P}(-2k-1)  (D)) \cong
R^0q_*(p^*\mathcal{O}_{\mathbb{P}}(k-2)(2D))^*$$
$$= S^{k-2}H^0(E)\otimes\mathcal{O}_S.$$
Similarly, we get
$$R^{k+2}q_* ( p^*{\mathcal
O}_{\mathbb
P}(-2k-1)\otimes\chi^*\mathcal{O}_{\mathbb{P}({\mathcal F})}(-1)
(D))$$
$$=R^{k+2}q_* ( p^*{\mathcal
O}_{\mathbb P}(-2k-2) (2D)) \cong
R^0q_*(p^*\mathcal{O}_{\mathbb{P}}(k-1)(D))^*$$
$$= S^{k-1}H^0(E)\otimes\mathcal{O}_S.$$
The exact sequence (\ref{le8ex3}) thus rewrites as :
$$ 0\rightarrow
{\mathcal V}\otimes L^{-1}(\Delta)\otimes S^{k-2}H^0(E)\rightarrow
R^{k+2}q_* (\chi^*{\mathcal H}\otimes p^*{\mathcal O}_{\mathbb
P}(-2k-1) (D)))$$ $$\rightarrow L^{-1}(\Delta)\otimes
S^{k-1}H^0(E)\rightarrow0.$$ This last exact sequence is now
canonically split because
$$H^0(S,\mathcal{V})=H^1(S,\mathcal{V})=0.$$
It follows that we have a canonical isomorphism :
\begin{eqnarray}\nonumber H^{k+4}(\widetilde{Y},\chi^*{\mathcal H}\otimes p^*{\mathcal
O}_{\mathbb P}(-2k-1) (D))=H^2(S,R^{k+2}q_* (\chi^*{\mathcal
H}\otimes p^*{\mathcal O}_{\mathbb P}(-2k-1) (D)))\\
 = H^2(S,{\mathcal V}\otimes L^{-1}(\Delta))\otimes
S^{k-2}H^0(E) \oplus H^2(S,L^{-1}(\Delta))\otimes S^{k-1}H^0(E)
\,\,\,\,\,.\label{point3}
\end{eqnarray}

In conclusion, using (\ref{point1}), (\ref{point2}) and
(\ref{point3}), we have found a canonical identification :
$$H^{k+3}(\widetilde{Y},S^{k+1}{\mathcal G}\otimes
q^*(L^{-1}(\Delta))(D))\cong$$ $$ Ker\,(S^{k-2}H^0(E)\otimes
H^2(S,\mathcal{V}\otimes L^{-1}(\Delta))\oplus
S^{k-1}H^0(E)\otimes H^2(S,L^{-1}(\Delta))$$
   $$\rightarrow S^{k-1}H^0(E)\otimes H^2(S,L^{-1}(\Delta))).$$
   One checks that the second component of the map is the
   identity, which gives a canonical identification :
$$H^{k+3}(\widetilde{Y},S^{k+1}{\mathcal G}\otimes
q^*(L^{-1}(\Delta))(D))\cong S^{k-2}H^0(E)\otimes
H^2(S,\mathcal{V}\otimes L^{-1}(\Delta)),$$ proving \ref{item1}.

The isomorphism \ref{item2} is proved in the same way.
 \cqfd

Let us now compute the spaces $H^2(S,\mathcal{V}\otimes
L^{-1}(\Delta))$ and $H^2(S,S^2\mathcal{V}\otimes
L^{-1}(\Delta))$.

  We first
observe that the exact sequence
$$0\rightarrow\mathcal{V}\rightarrow H^0(E(-\Delta))\otimes
\mathcal{O}_S\rightarrow E(-\Delta)\rightarrow0,
$$
and Serre duality give an isomorphism
\begin{eqnarray}\label{isoNu}H^2(S,\mathcal{V}\otimes L^{-1}(\Delta))\cong Ker\,c,
\end{eqnarray} where $c$ is the contraction map
$$H^0(E(-\Delta))\otimes H^0(S,L(-\Delta))^*\rightarrow
H^0(E(-\Delta))^*$$ induced by the determinant map
$$d_0:\bigwedge^2H^0(E(-\Delta))\rightarrow H^0(S,L(-\Delta)).$$
Similarly the induced exact sequence
$$
0\rightarrow S^2\mathcal{V}\rightarrow S^2H^0(E(-\Delta))\otimes
\mathcal{O}_S\rightarrow H^0(E(-\Delta))\otimes E(-\Delta)
\rightarrow\bigwedge^2E(-\Delta)\rightarrow0 $$ gives a surjective
map :
 \begin{eqnarray}\label{isoNu2}H^2(S,S^2{\mathcal V}\otimes
L^{-1}(\Delta))\twoheadrightarrow Ker\,c',
\end{eqnarray}
where the contraction map
$$c':S^2H^0(E(-\Delta))\otimes H^0(S,L(-\Delta))^*\rightarrow
H^0(E(-\Delta))\otimes H^0(E(-\Delta))^*,$$  is also induced by
$d_0$.

Using Lemma \ref{le1mardi29}, we   rewrite now the sequence
(\ref{complexe9mai}) as follows:
 We first identify
$$\bigwedge^{2k+1}H^0(L)^*,\,
\bigwedge^{2k}H^0(L)^*,\, \bigwedge^{2k-1}H^0(L)^*$$ respectively to
$$H^0(L),\,\bigwedge^2H^0(L),\, \bigwedge^3H^0(L).$$
 Then via the isomorphisms given in Lemma \ref{le1mardi29}, and
(\ref{isoNu}), (\ref{isoNu2})
  above, our sequence (\ref{complexe9mai}) becomes,
 after replacing the first term by its quotient given in (\ref{isoNu2}),
 through which the first map factors :
 \begin{eqnarray}\label{K1}
 H^0(L)\otimes S^{k-1}H^0(E)\otimes Ker\,c'
 \rightarrow \bigwedge^2H^0(L)\otimes S^{k-2}H^0(E)\otimes
 Ker\,c\\
 \nonumber
 \rightarrow \bigwedge^3H^0(L)\otimes S^{k-3}H^0(E)
 \otimes H^0(S,L(-\Delta))^*.
 \end{eqnarray}
 It is immediate to check that the maps of the complex are
 induced by
the determinant map
$$d:H^0(E)\otimes H^0(E(-\Delta))\rightarrow H^0(L)$$
and by the natural maps, for $(i,j)=(k-1,2),\,(k-2,1),\,(k-3,0)$ :
$$ S^iH^0(E)\otimes S^jH^0(E(-\Delta)) \rightarrow
S^{i-1}H^0(E)\otimes S^{j-1}H^0(E(-\Delta)) \otimes H^0(E)\otimes
H^0(E(-\Delta)).$$

\vspace{0,5cm}

{\it Fourth step.} We now want to prove Lemma \ref{le2mardi29},
which we have just proved to be equivalent to exactness at the
middle of the sequence (\ref{K1}). We do not need at this point
the $K3$ surface anymore. We shall use now only the map $d$ and do
geometry on the Grassmannian of subspaces of $H^0(E)$. We believe
that this step is the only essential one in the proof of
proposition \ref{injective}.

   Denote by $G$ the Grassmannian of rank $2$ subspaces of
$H^0(E)$, and let $\mathcal{L}$ be the Pl\"ucker line bundle on
$G$,  $\mathcal{E}$ the tautological rank $2$ quotient bundle on
$G$. Let $ \widetilde{G}'$ be the desingularization of the
hypersurface $G'\subset G$ parametrizing the $V\subset H^0(E)$
meeting $H^0(E(-\Delta))$, defined as $$ \widetilde{G}'=\{
(v,V)\in \mathbb{P}(H^0(E(-\Delta)))\times G,\, v\in V\}.
$$ We shall also denote by $\mathcal{L},\, \mathcal{E}$ the
pull-backs of $\mathcal{L},\,\mathcal{E}$ to $\widetilde{G}'$ by
the second projection.
 Let
$$g:\widetilde{G}'\rightarrow \mathbb{P}H^0(E(-\Delta))$$
be the first projection, and denote by $H$ the line bundle
$g^*\mathcal{O}_{\mathbb{P}H^0(E(-\Delta))}(1)$.

Next, the map
$$d_0:\bigwedge^2H^0(E(-\Delta))\rightarrow H^0(L(-\Delta))$$
allows to define
$$I\subset
\mathbb{P}H^0(L(-\Delta))^*\times\mathbb{P}H^0(E(-\Delta)),$$
$$I=\{(\sigma,v),\,\sigma(d_0(v\wedge\,\cdot))=0\,\,{\rm in}\,\,H^0(E(-\Delta))^*\}.$$
Let $\widetilde{I}\stackrel{\pi}{\rightarrow}\widetilde{G}'$ be
the fibered product
$I\times_{\mathbb{P}H^0(E(-\Delta))}\widetilde{G}'$. We shall
denote by $\mathcal{K}$ the line bundle
$$pr_1^*\mathcal{O}_{\mathbb{P}H^0(L(-\Delta))^*}(1)$$ on
$\widetilde{I}$, and by $\mathcal{L}',\, \mathcal{E}', H'$ the
pull-backs to $\widetilde{I}$ of the corresponding objects on
$\widetilde{G}'$, that is
$$\mathcal{L}'=\pi^*\mathcal{L},\,\,\,\mathcal{E}'=\pi^*\mathcal{E},\,\,\,\,H'=\pi^*H.$$
 The notations are summarized in the
following diagram:
\begin{eqnarray}\label{diagram2}
\begin{diagram}
\node{I}\arrow{s,l}{pr_1}\arrow{se,l}{pr_2}\node{\widetilde{I}}\arrow{w,t}{g'}\arrow{se,l}{\pi}
\node{P_{\widetilde{I}}}
\arrow{w,t}{\alpha}\arrow{se}\arrow{see,l}{\beta}\node{}\\
\node{\mathbb{P}H^0(L(-\Delta))^*}\node{\mathbb{P}H^0(E(-\Delta))}
\node{\widetilde{G}'}\arrow{w,t}{g}\arrow{s,l}{desing}
\arrow{se}\node{P}\arrow{s,l}{p}\arrow{e,b}{q}\node{\mathbb{P}H^0(E)}\\
\node{}\node{}\node{G'}\arrow{e,b,J}{hypersurface}\node{G}\node{}
\end{diagram}
\end{eqnarray}

We shall use the following Lemma:
\begin{lemm} \label{utileplustardaussi}
For any positive integers $p,\,s,\,t\geq s,\,p\geq s$ we have
$$ H^0(\widetilde{I},S^{p}\mathcal{E}'\otimes\mathcal{L}'^{-s}\otimes
H'^t\otimes \mathcal{K})\cong
H^0(\widetilde{I},S^{p-s}\mathcal{E}'\otimes H'^{t-s}\otimes
\mathcal{K}),$$
$$H^0(\widetilde{I},S^{p-s}\mathcal{E}'\otimes
H'^{t-s}\otimes \mathcal{K})\cong S^{p-s}H^0(E)^*\otimes
[S^{t-s}H^0(E(-\Delta))\otimes H^0(L(-\Delta))^*]_0^*.$$
\end{lemm}
Here, the term  $[S^{t-s}H^0(E(-\Delta))\otimes
H^0(L(-\Delta))^*]_0$ has the following meaning: the map $d_0$
provides a contraction map $$H^0(E(-\Delta))\otimes
H^0(L(-\Delta))^*\rightarrow  H^0(E(-\Delta))^*$$ and more
generally a contraction map
$$ c_i: S^iH^0(E(-\Delta))\otimes
H^0(L(-\Delta))^*\rightarrow S^{i-1}H^0(E(-\Delta))\otimes
H^0(E(-\Delta))^*.$$ Then we denote
$$[S^{i}H^0(E(-\Delta))\otimes
H^0(L(-\Delta))^*]_0:=Ker\,c_{i}.$$ Note that with the previous
notations, we have $c'=c_2,\,c=c_1$.

\vspace{0,5cm}

 {\bf Proof of Lemma \ref{utileplustardaussi}.}
 These facts are proved using the exact sequence on
 $\widetilde{G}'$:
 $$0\rightarrow \mathcal{L}\otimes H^{-1}\rightarrow
 \mathcal{E}\rightarrow H\rightarrow0,$$
 or, pulling-back via $\pi$:
 $$0\rightarrow \mathcal{L}'\otimes H'^{-1}\rightarrow
 \mathcal{E}'\rightarrow H'\rightarrow0$$
 on $\widetilde{I}$.
 Taking symmetric powers
 gives
    exact sequences
 \begin{eqnarray}\label{syml10june}
 0\rightarrow S^{l-1}\mathcal{E}'\otimes\mathcal{L}'\otimes H'^{-1}\rightarrow
 S^l\mathcal{E}'\rightarrow H'^l\rightarrow0
 \end{eqnarray}
 on $\widetilde{I}$.
 Take $l=p$ and tensor (\ref{syml10june}) with
 $\mathcal{L}'^{-s}\otimes H'^t\otimes \mathcal{K}$. Observing that if $s>0$,
 $ H^0(\widetilde{I},\mathcal{L}'^{-s}\otimes H'^{p+t}\otimes \mathcal{K})=0
 $ because
 the restriction of this line bundle to the fibers of $g'$ is
 $\mathcal{O}(-s)$ on a projective space of dimension $>0$, we
 conclude that
 $$H^0(\widetilde{I},S^{p}\mathcal{E}'\otimes\mathcal{L}'^{-s}\otimes
H'^t\otimes \mathcal{K})\cong
H^0(\widetilde{I},S^{p-1}\mathcal{E}'\otimes\mathcal{L}'^{-s+1}\otimes
H'^{t-1}\otimes \mathcal{K}),$$ which proves the first equality by
iteration.

  The second equality   follows
 from the following observation : denoting by
 $P_{\widetilde{I}}\stackrel{(\alpha,\beta)}{\rightarrow}
  \widetilde{I}\times \mathbb{P}(H^0(E))$  the pull-back to $\widetilde{I}$
  of the tautological
 $\mathbb{P}^1$-bundle $P$ on $G$, (see diagram (\ref{diagram2}),)
 there is a natural  map
 $$
  P_{\widetilde{I}}\stackrel{(g'\circ\alpha,\beta)}{\longrightarrow}
   I\times\mathbb{P}(H^0(E))$$
   which is immediately seen to be birational.
 Furthermore, we have $S^l\mathcal{E}'\cong
 R^0\alpha_*(\beta^*\mathcal{O}(l))$ on $\widetilde{I} $, so that
 $$H^0(\widetilde{I},H'^i\otimes S^j\mathcal{E}'\otimes \mathcal{K})=
 H^0(P_{\widetilde{I}},\alpha^*(H'^i\otimes
 \mathcal{K})\otimes\beta^*\mathcal{O}(j))$$
 $$=H^0(P_{\widetilde{I}},(g'\circ\alpha)^*(pr_2^*\mathcal{O}(i)\otimes \mathcal{K}))
 \otimes\beta^*\mathcal{O}(j)).$$
 Because the map $(g'\circ\alpha,\beta) $ is birational, this is
 also equal to
 $$ H^0(I\times\mathbb{P}(H^0(E)),pr_2^*\mathcal{O}(i)\otimes \mathcal{K}\boxtimes
 \mathcal{O}(j))=
  H^0(I,pr_2^*\mathcal{O}(i)\otimes \mathcal{K})\otimes
  S^jH^0(E)^*.$$
To conclude, it thus  suffices to show the following equality:
 $$H^0(I,pr_2^*\mathcal{O}(i)\otimes \mathcal{K})
 =[S^iH^0(E(-\Delta))\otimes H^0(L(-\Delta))^*]_0^*.$$
 This last fact follows from the fact that by definition of $I$
 and $\mathcal{K}$,
 the vector bundle $$\mathcal{Q}:=R^0pr_{2*}\mathcal{K}$$
 on $\mathbb{P}H^0(E(-\Delta))$ fits into the exact sequence
 \begin{eqnarray}\label{Q4aout}
 0\rightarrow \mathcal{O}(-2)\rightarrow \mathcal{O}(-1)\otimes
 H^0(E(-\Delta))\stackrel{d_0}{\rightarrow}
 H^0(L(-\Delta))\otimes\mathcal{O}\rightarrow\mathcal{Q}
 \rightarrow 0.
 \end{eqnarray}
It follows  from this exact sequence, using the fact that
$rk\,H^0(E(-\Delta))=k+1\geq3$ and vanishing on
$\mathbb{P}H^0(E(-\Delta))$, that
$$H^0(I,pr_2^*\mathcal{O}(i)\otimes
\mathcal{K})=H^0(\mathbb{P}H^0(E(-\Delta)),\mathcal{O}(i)\otimes
R^0pr_{2*}\mathcal{K})$$
$$=H^0(\mathbb{P}H^0(E(-\Delta)),\mathcal{O}(i)\otimes
\mathcal{Q})$$ is equal to the cokernel of the map induced by
$d_0$
$$S^{i-1}H^0(E(-\Delta))^*\otimes H^0(E(-\Delta))\rightarrow
S^{i}H^0(E(-\Delta))^*\otimes H^0(L(-\Delta)),$$ that is to
$[S^iH^0(E(-\Delta))\otimes H^0(L(-\Delta))^*]_0^*$.

 \cqfd This Lemma will provide in particular canonical
identifications:

\begin{eqnarray}\label{id1aout}
H^0(\widetilde{I},S^{k-2}\mathcal{E}'\otimes\mathcal{K}\otimes
H'\otimes \mathcal{L}'^{-1}) =S^{k-3}H^0(E)^*\otimes
H^0(L(-\Delta)),
\end{eqnarray}
\begin{eqnarray}\label{id2aout}
H^0(\widetilde{I},S^{k-2}\mathcal{E}'\otimes\mathcal{K}\otimes
H')=\\ \nonumber S^{k-2}H^0(E)^*\otimes [H^0(E(-\Delta))\otimes
H^0(L(-\Delta))^*]_0^* \\ \nonumber= S^{k-2}H^0(E)^*\otimes
(Ker\,c)^*,\end{eqnarray} an inclusion
\begin{eqnarray}
\label{inclusion3aout}H^0(\widetilde{I},S^{k-2}\mathcal{E}'\otimes\mathcal{K}
\otimes H'\otimes\mathcal{L}')\subset
H^0(\widetilde{I},S^{k-1}\mathcal{E}'\otimes\mathcal{K} \otimes
H'^2)
\end{eqnarray}
and an identification:
\begin{eqnarray}\label{id3aout}
H^0(\widetilde{I},S^{k-1}\mathcal{E}'\otimes\mathcal{K} \otimes
H'^2)=\\ \nonumber S^{k-1}H^0(E)^*\otimes
[S^2H^0(E(-\Delta))\otimes H^0(L(-\Delta))^*]_0^* \\ \nonumber=
S^{k-1}H^0(E)^*\otimes (Ker\,c')^*.
\end{eqnarray} They are used as follows: The determinant map
$d:H^0(E)\otimes H^0(E(-\Delta))\rightarrow H^0(L)$ provides
dually a linear system
$$H^0(L)^*=: W\subset H^0(G',\mathcal{L})\subset
H^0(\widetilde{I},\mathcal{L}')$$ which has no base-point by
Proposition \ref{detgrass}. This provides an exact Koszul complex
on $\widetilde{I}$ :
\begin{eqnarray} \label{K13mai}0\rightarrow
\bigwedge^{2k+2}W\otimes\mathcal{L}'^{-2k-2}\rightarrow\ldots
\rightarrow W\otimes\mathcal{L}'^{-1}\rightarrow
\mathcal{O}_{\widetilde{I}}\rightarrow0. \end{eqnarray}
 We twist
by $S^{k-2}\mathcal{E}'\otimes \mathcal{L}'^{\otimes2} \otimes
H'\otimes\mathcal{K}$ and take global sections.
 The relevant piece of this   complex of global sections is:
 \begin{eqnarray}
 \label{Koszul9mai}\ldots\rightarrow \bigwedge^3W\otimes
 H^0(\widetilde{I},S^{k-2}\mathcal{E}'\otimes\mathcal{K}\otimes
H'\otimes \mathcal{L}'^{-1}) \rightarrow
\\
\nonumber\bigwedge^2W\otimes
H^0(\widetilde{I},S^{k-2}\mathcal{E}'\otimes\mathcal{K}\otimes
H')\rightarrow W\otimes
H^0(\widetilde{I},S^{k-2}\mathcal{E}'\otimes\mathcal{K}\otimes
H'\otimes \mathcal{L}')\rightarrow\ldots
\end{eqnarray}
Using the inclusion (\ref{inclusion3aout}), this sequence has the
same cohomology at the middle as the sequence:
\begin{eqnarray}
 \label{Koszul9maimodif}\ldots\rightarrow \bigwedge^3W\otimes
 H^0(\widetilde{I},S^{k-2}\mathcal{E}'\otimes\mathcal{K}\otimes
H'\otimes \mathcal{L}'^{-1}) \rightarrow
\\
\nonumber\bigwedge^2W\otimes
H^0(\widetilde{I},S^{k-2}\mathcal{E}'\otimes\mathcal{K}\otimes
H')\rightarrow W\otimes
H^0(\widetilde{I},S^{k-1}\mathcal{E}'\otimes\mathcal{K}\otimes
H'^2)\rightarrow\ldots
\end{eqnarray}
 Finally, using the
identifications (\ref{id1aout}), (\ref{id2aout}), (\ref{id3aout})
above, we see that the three terms of this last sequence are
canonically dual to the three terms of the sequence (\ref{K1}). We
leave to the reader to verify that this last sequence is indeed
dual to (\ref{K1}).
 Hence, the exactness at the middle of
 (\ref{K1}) is equivalent to the exactness at the middle of the
 Koszul
sequence (\ref{Koszul9mai}), and we claim that this is implied by
the following  statement :

\vspace{0,5cm}

\begin{lemm}\label{derlemaout} For $1\leq i\leq 2k-1$ we have the vanishing
\begin{eqnarray}\label{penible12juin}
H^i(\widetilde{I},S^{k-2}\mathcal{E}'\otimes \mathcal{K}\otimes
H'\otimes\mathcal{L}'^{-i-1})=0.
\end{eqnarray}
\end{lemm}

 Indeed, the Koszul complex (\ref{K13mai}) twisted
by $ S^{k-2}\mathcal{E}'\otimes\mathcal{K} \otimes H'
\otimes\mathcal{L}'^2   $ reads
$$0\rightarrow
\bigwedge^{2k+2}W\otimes S^{k-2}\mathcal{E}'\otimes \mathcal{K}
\otimes H'\otimes\mathcal{L}'^{-2k}\rightarrow\ldots \rightarrow
\bigwedge^2W\otimes S^{k-2}\mathcal{E}'\otimes\mathcal{K} \otimes
H'$$
$$\rightarrow W\otimes S^{k-2}\mathcal{E}'\otimes\mathcal{K}\otimes
 H'\otimes \mathcal{L}' \rightarrow\ldots.
$$
In  particular, this provides by truncation a resolution
$$0\rightarrow \bigwedge^{2k+2}W\otimes S^{k-2}\mathcal{E}'\otimes
\mathcal{K} \otimes H'\otimes\mathcal{L}'^{-2k}\rightarrow\ldots
\rightarrow \bigwedge^4W\otimes
S^{k-2}\mathcal{E}'\otimes\mathcal{K} \otimes H'\otimes
\mathcal{L}'^{-2} \rightarrow \mathcal{M}\rightarrow 0,$$ of the
sheaf $\mathcal{M}$ on $\widetilde{I}$ which fits in the exact
sequence:
\begin{eqnarray}0\rightarrow\mathcal{M}\rightarrow
\bigwedge^3W\otimes S^{k-2}\mathcal{E}'\otimes\mathcal{K} \otimes
H'\otimes\mathcal{L}'^{-1}\\ \nonumber\rightarrow
\bigwedge^2W\otimes S^{k-2}\mathcal{E}'\otimes\mathcal{K} \otimes
H' \rightarrow W\otimes S^{k-2}\mathcal{E}'\otimes\mathcal{K}
\otimes H'\otimes \mathcal{L}'\ldots.
\end{eqnarray}
Now the exactness at the middle of the sequence (\ref{Koszul9mai})
is implied by the vanishing $H^1(\widetilde{I},\mathcal{M})=0$. By
the above resolution, this space is isomorphic to the
hypercohomology group
\begin{eqnarray}\label{hypercoh12juin}\mathbb{H}^1(\widetilde{I},0\rightarrow
\bigwedge^{2k+2}W\otimes S^{k-2}\mathcal{E}'\otimes \mathcal{K}
\otimes H'\otimes\mathcal{L}'^{-2k}\rightarrow\ldots \\
\nonumber\rightarrow \bigwedge^4W\otimes
S^{k-2}\mathcal{E}'\otimes\mathcal{K} \otimes H'\otimes
\mathcal{L}'^{-2} \rightarrow0), \end{eqnarray} where the last
term on the right is put in degree $0$. The terms $E_1^{p,q}$ of
the spectral sequence associated to the naive filtration of this
complex are equal to
$$ \bigwedge^{4-p}W\otimes
H^q(\widetilde{I},S^{k-2}\mathcal{E}'\otimes\mathcal{K}\otimes
H'\otimes\mathcal{L}'^{p-2})$$ in degree $p+q$, that is, for
$p+q=1$, to
$$\bigwedge^{q+3}W\otimes
H^q(\widetilde{I},S^{k-2}\mathcal{E}'\otimes\mathcal{K}\otimes
H'\otimes\mathcal{L}'^{-q-1}),\,q\geq1.$$
 Hence the vanishings
(\ref{penible12juin}) say that the terms $E_1^{p,q},\,p+q=1$ of
this spectral sequence  are $0$ for $1\leq q\leq 2k-1$ and they
are also obviously $0$ for $q>2k-1$, since then
$\bigwedge^{q+3}W=0$. Thus (\ref{hypercoh12juin}) vanishes and so
does $H^1(\widetilde{I},\mathcal{M})$. \cqfd {\bf Proof of Lemma
\ref{derlemaout}.}  Consider the Cartesian diagram (see diagram
(\ref{diagram2}):
\begin{eqnarray}\begin{diagram}
\node{}\node{\widetilde{I}}\arrow{e,t}{\pi}\arrow{s}\node{\widetilde{G}'}
\arrow{s,l}{g}
\\\node{\mathbb{P}H^0(L(-\Delta))^*} \node{I}\arrow{w,b}{pr_1}\arrow{e,b}{pr_2}
\node{\mathbb{P}H^0(E(-\Delta))}
\end{diagram}
\end{eqnarray}
We have
$$R^0\pi_*\mathcal{K}=g^*(R^0pr_{2*}(pr_1^*\mathcal{O}(1)))=g^*\mathcal{Q},$$
where the bundle $\mathcal{Q}$ on $\mathbb{P}H^0(E(-\Delta))$
admits the resolution (see (\ref{Q4aout})):
$$
 0\rightarrow \mathcal{O}(-2)\rightarrow \mathcal{O}(-1)\otimes
 H^0(E(-\Delta))\stackrel{d_0}{\rightarrow}
 H^0(L(-\Delta))\otimes\mathcal{O}\rightarrow\mathcal{Q}
 \rightarrow 0.
 $$
 Since $\mathcal{E}',\,\mathcal{L}',\,H'$ are pull-backs via $\pi$ of the
 corresponding objects on $\widetilde{G}'$, we have:
 $$H^i(\widetilde{I},S^{k-2}\mathcal{E}'\otimes \mathcal{K}\otimes
H'\otimes\mathcal{L}'^{-i-1})=
H^i(\widetilde{G}',S^{k-2}\mathcal{E}\otimes
H\otimes\mathcal{L}^{-i-1}\otimes g^*\mathcal{Q}).$$ The bundle
$g^*\mathcal{Q}$ admits the resolution:
$$ 0\rightarrow H^{-2}\rightarrow H^{-1}\otimes
 H^0(E(-\Delta))\stackrel{d_0}{\rightarrow}
 H^0(L(-\Delta))\otimes\mathcal{O}_{\widetilde{G}'}\rightarrow g^*\mathcal{Q}
 \rightarrow 0,
 $$
 and it follows that the desired  vanishing
 $$H^i(\widetilde{G}',S^{k-2}\mathcal{E}\otimes
H\otimes\mathcal{L}^{-i-1}\otimes g^*\mathcal{Q})=0,\,\, 1\leq
i\leq 2k-1$$ is a consequence of the following:
\begin{eqnarray}\label{van34aout}H^i(\widetilde{G}',S^{k-2}\mathcal{E}\otimes
H\otimes\mathcal{L}^{-i-1})=0,\,\, 1\leq i\leq 2k-1,\\ \nonumber
H^{i+1}(\widetilde{G}',S^{k-2}\mathcal{E}\otimes\mathcal{L}^{-i-1})=0,\,\,
1\leq i\leq 2k-1,
\\ \nonumber
H^{i+2}(\widetilde{G}',S^{k-2}\mathcal{E}\otimes
H^{-1}\otimes\mathcal{L}^{-i-1})=0,\,\, 1\leq i\leq 2k-2,
\end{eqnarray}
and of the following fact:

\begin{lemm} \label{derder4aout}The  map
$$H^{i+2}(\widetilde{G}',S^{k-2}\mathcal{E}\otimes
H^{-1}\otimes\mathcal{L}^{-i-1})\rightarrow H^0(E(-\Delta))\otimes
H^{i+2}(\widetilde{G}',S^{k-2}\mathcal{E}\otimes\mathcal{L}^{-i-1})$$
induced by the natural inclusion $$H^{-1}\subset
H^0(E(-\Delta))\otimes\mathcal{O}_{\widetilde{G}'}$$ is injective
for $i=2k-1$.
\end{lemm}

Let us first prove (\ref{van34aout}): We use  the fact that we can
see $\widetilde{G}'$ as the complete intersection of two members
of $\mid H\mid$ on the tautological $\mathbb{P}^1$-bundle $P$ on
$G$. Indeed, let $P\subset G\times \mathbb{P}H^0(E)$ be the
tautological subbundle, and denote by $p:P\rightarrow G$ the first
projection, $q:P\rightarrow \mathbb{P}H^0(E)$ the second
projection (see diagram (\ref{diagram2})). Denote by $H$ the line
bundle $q^*\mathcal{O}(1)$ on $P$, and by
$\mathcal{E},\,\mathcal{L}$ the pull-backs via $p$ of the
corresponding bundles on $G$. Then by definition, $\widetilde{G}'$
identifies to $q^{-1}(\mathbb{P}H^0(E(-\Delta)))$, and the bundles
$H,\,\mathcal{E},\,\mathcal{L}$ are the restrictions to
$\widetilde{G}'$ of the corresponding objects on $P$.

Hence there is a Koszul resolution of $\mathcal{O}_{\widetilde{G}'}$
which has the form:
$$0\rightarrow\bigwedge^2R\otimes H^{-2}\rightarrow R\otimes H^{-1}
\rightarrow\mathcal{O}_{P}\rightarrow\mathcal{O}_{\widetilde{G}'}
\rightarrow0,$$ where $R$ is a rank $2$ vector space.

Using this resolution, we see that the vanishing statements
(\ref{van34aout}) are a consequence of the following ones:
\begin{enumerate}
\item\label{1vanaout}$H^i(P,S^{k-2}\mathcal{E}\otimes
H\otimes\mathcal{L}^{-i-1})=0,\,\, 1\leq i\leq 2k-1$,
\item\label{2vanaout}$H^{i+1}(P,S^{k-2}\mathcal{E}\otimes\mathcal{L}^{-i-1})=0,\,\,
1\leq i\leq 2k-1$,
\item\label{3vanaout}$H^{i+2}(P,S^{k-2}\mathcal{E}\otimes
H^{-1}\otimes\mathcal{L}^{-i-1})=0,\,\, 1\leq i\leq 2k-1$,
\item\label{4vanaout}$H^{i+3}(P,S^{k-2}\mathcal{E}\otimes
H^{-2}\otimes\mathcal{L}^{-i-1})=0,\,\, 1\leq i\leq 2k-1$,
\item\label{5vanaout}$H^{i+4}(P,S^{k-2}\mathcal{E}\otimes
H^{-3}\otimes\mathcal{L}^{-i-1})=0,\,\, 1\leq i\leq 2k-2$.
\end{enumerate}
Recall now the following statement proven in the Appendix of
\cite{Vo}:
\begin{prop} \label{propappen}For $q>0,\,q'\geq0$, $H^p(G,\mathcal{L}^{-q}\otimes
S^{q'}\mathcal{E})=0$ if $p\not=k+1,\,2k+2$.

For $p=k+1$, $H^p(G,\mathcal{L}^{-q}\otimes S^{q'}\mathcal{E})=0$ if
$-q+q'+1<0$ or $q\leq k+1$.
\end{prop}
(Note the shift of notation from $k$ there to $k+1$ here, which is
due to the fact that we are now working with a space $H^0(E)$ of
rank $k+3$ instead of $k+2$.)

The vanishing \ref{2vanaout} follows directly from this
Proposition. The vanishing \ref{3vanaout} follows from the fact
that $H^{-1}$ has trivial cohomology along the fibers of
$p:P\rightarrow G$. For the vanishing \ref{1vanaout}, we use  the
exact sequence on $P$:
$$0\rightarrow\mathcal{L}\otimes H^{-1}\rightarrow
\mathcal{E}\rightarrow H\rightarrow0.$$ It provides the exact
sequence:
$$0\rightarrow\mathcal{L}\otimes H^{-1}\otimes S^{k-3}\mathcal{E}\rightarrow
S^{k-2}\mathcal{E}\rightarrow H^{k-2}\rightarrow0.$$ Hence we see
that \ref{1vanaout} is implied by the vanishings:
$$H^i(P,S^{k-3}\mathcal{E}\otimes\mathcal{L}^{-i})=0,\,\,1\leq
i\leq2k-1,$$
$$H^i(P,H^{k-1}\otimes\mathcal{L}^{-i-1})=H^i(G,S^{k-1}\mathcal{E}\otimes\mathcal{L}^{-i-1})=0,\,\,1\leq
i\leq2k-1,$$ and they are both consequences of Proposition
\ref{propappen}.

For the vanishing \ref{4vanaout}, one notes that $K_{P/G}$ is equal
to $H^{-2}\otimes\mathcal{L}$. Hence we have
$$H^{i+3}(P,S^{k-2}\mathcal{E}\otimes\mathcal{L}^{-i-1}\otimes
H^{-2})=
H^{i+2}(G,S^{k-2}\mathcal{E}\otimes\mathcal{L}^{-i-1}\otimes
R^1p_*H^{-2})$$
$$=H^{i+2}(G,S^{k-2}\mathcal{E}\otimes\mathcal{L}^{-i-2}).$$
This vanishes for $1\leq i\leq2k-1$ by Proposition
\ref{propappen}.

To conclude, \ref{5vanaout} is proved as follows: We have as above:
$$H^{i+4}(P,S^{k-2}\mathcal{E}\otimes\mathcal{L}^{-i-1}\otimes
H^{-3})=
H^{i+3}(G,S^{k-2}\mathcal{E}\otimes\mathcal{L}^{-i-1}\otimes
R^1p_*H^{-3})$$
$$=H^{i+3}(G,S^{k-2}\mathcal{E}\otimes\mathcal{L}^{-i-2}\otimes\mathcal{E}^*),$$
using relative Serre duality and $R^0p_*H=\mathcal{E}$ on $G$.
Since $\mathcal{E}^*=\mathcal{E}\otimes\mathcal{L}^{-1}$, the last
term is equal to
$$H^{i+3}(G,S^{k-2}\mathcal{E}\otimes\mathcal{E}\otimes\mathcal{L}^{-i-3}).$$
By the exact sequence
$$0\rightarrow\mathcal{L}\otimes S^{k-3}\mathcal{E}\rightarrow
S^{k-2}\mathcal{E}\otimes\mathcal{E}\rightarrow S^{k-1}\mathcal{E}
\rightarrow0,$$ we see that  the vanishing
$H^{i+3}(G,S^{k-2}\mathcal{E}\otimes\mathcal{E}\otimes\mathcal{L}^{-i-3})=0$
for $1\leq i\leq2k-2$  is a consequence of the vanishings
$$H^{i+3}(G,S^{k-3}\mathcal{E}\otimes\mathcal{L}^{-i-2})=0,$$
$$H^{i+3}(G,S^{k-1}\mathcal{E}\otimes\mathcal{L}^{-i-3})=0,$$
for $1\leq i\leq2k-2$, which both follow from Proposition
\ref{propappen}.

 This concludes the proof of  (\ref{van34aout}) and
 the proof of Lemma \ref{derlemaout} will then be concluded with the
 proof of Lemma \ref{derder4aout}.
 \cqfd
{\bf Proof of Lemma \ref{derder4aout}.}

 Since the map
$$H^{2k+1}(\widetilde{G}',S^{k-2}\mathcal{E}\otimes
H^{-1}\otimes\mathcal{L}^{-2k})\rightarrow H^0(E(-\Delta))\otimes
H^{2k+1}(\widetilde{G}',S^{k-2}\mathcal{E}\otimes\mathcal{L}^{-2k})$$
is induced by the inclusion
$$H^{-1}\subset H^0(E(-\Delta))\otimes\mathcal{O}_{\widetilde{G}'},$$
which is dual to the evaluation map, where $H^0(E(-\Delta))^*$ is
identified to $H^0(\widetilde{G}',H)$, we see that its dual is
equal to the multiplication map:
$$ H^0(E(-\Delta))^*\otimes H^0(\widetilde{G}',K_{\widetilde{G}'}\otimes
S^{k-2}\mathcal{E}^*\otimes\mathcal{L}^{2k}) \rightarrow
H^0(\widetilde{G}',K_{\widetilde{G}'}\otimes H\otimes
S^{k-2}\mathcal{E}^*\otimes\mathcal{L}^{2k}).$$

The canonical bundle of $\widetilde{G}'$ is equal to
$\mathcal{L}^{-k-2}$, because $\widetilde{G}'$ is the complete
intersection of two members of $\mid H\mid$ in $P$ and
$K_P=\mathcal{L}^{-k-2}\otimes H^{-2}$.

Thus we have to prove that the multiplication map
$$H^0(E(-\Delta))^*\otimes H^0(\widetilde{G}',
S^{k-2}\mathcal{E}^*\otimes\mathcal{L}^{k-2}) \rightarrow
H^0(\widetilde{G}', H\otimes
S^{k-2}\mathcal{E}^*\otimes\mathcal{L}^{k-2})$$ is surjective.

Since $\mathcal{E}^*\otimes\mathcal{L}\cong \mathcal{E}$, this is
equivalent to the surjectivity of the multiplication map
$$
H^0(E(-\Delta))^*\otimes H^0(\widetilde{G}', S^{k-2}\mathcal{E})
\rightarrow H^0(\widetilde{G}', H\otimes S^{k-2}\mathcal{E}).
$$
This follows from the surjectivity of the multiplication map
$$H^0(P,H)\otimes H^0(P,S^{k-2}\mathcal{E})\rightarrow
H^0(P,H\otimes S^{k-2}\mathcal{E}),$$ and  of the restriction map
$$H^0(P,H\otimes S^{k-2}\mathcal{E})\rightarrow
H^0(\widetilde{G}', H\otimes S^{k-2}\mathcal{E}).$$
 \cqfd

\end{document}